\documentclass[11pt]{article}

\usepackage{latexsym}
\usepackage{amssymb}
\usepackage{amsfonts}
\usepackage{amsmath}
\usepackage{amsthm}
\usepackage{mathrsfs}
\usepackage{dsfont}    
\usepackage{bbold}     
\usepackage[english]{babel}
\usepackage{caption}
\usepackage{epsfig}
\usepackage{float}
\usepackage{subfigure}
\usepackage{psfrag}
\usepackage{graphicx}
\usepackage{epsfig}

\usepackage{ifthen}
\newboolean{pdflat}
\setboolean{pdflat}{true}

\newtheorem*{prop*}{Theorem}

\newcommand{\zerarcounters}{\setcounter{equation}{0}\setcounter{theorem}{0}\setcounter{figure}{0}\setcounter{table}{0}}

\newcommand{\st}{\scriptstyle}
\newcommand{\ZZZ}{\mathds{Z}}

\newcommand{\NNN}{\mathds{N}}
\newcommand{\QQQ}{\mathds{Q}}
\newcommand{\RRR}{\mathds{R}}
\newcommand{\TTT}{\mathds{T}}
\newcommand{\uno}{\mathds{1}}

\newcommand{\calQ}{{\mathcal Q}}

\newcommand{\e}{\varepsilon}
\newcommand{\al}{\alpha}
\newcommand{\de}{\delta}

\newcommand{\z}{\zeta}

\newcommand{\m}{\mu}

\newcommand{\g}{\gamma}
\newcommand{\om}{\omega}
\newcommand{\h}{\eta}
\newcommand{\la}{\lambda}
\newcommand{\f}{\varphi}
\newcommand{\s}{\sigma}

\newcommand{\der}{{\rm d}}
\newcommand{\ii}{{\rm i}}

\newcommand{\sn}{{\rm sn\,}}
\newcommand{\cn}{{\rm cn\,}}
\newcommand{\dn}{{\rm dn\,}}

\def\ins#1#2#3{\vbox to0pt{\kern-#2 \hbox{\kern#1 #3}\vss}\nointerlineskip}

\begin{document}


\title{\bf Attractiveness of periodic orbits in parametrically forced systems
with time-increasing friction}
\author
{\bf Michele Bartuccelli$^1$, Jonathan Deane$^1$, Guido Gentile$^2$
\vspace{2mm}
\\ \small 
$^1$ Department of Mathematics, University of Surrey, Guildford, GU2 7XH, UK 
\\ \small
$^2$ Dipartimento di Matematica, Universit\`a di Roma Tre, Roma, I-00146, Italy
\\ \small 
E-mail: m.bartuccelli@surrey.ac.uk, j.deane@surrey.ac.uk, gentile@mat.uniroma3.it
}


\date{}

\maketitle

\begin{abstract}
We consider dissipative one-dimensional systems subject to a periodic force
and study numerically how a time-varying friction affects the dynamics.
As a model system, particularly suited for numerical analysis,
we investigate the driven cubic oscillator in the presence of friction. 
We find that, if the damping coefficient increases in time up to a final
constant value, then the basins of attraction of the
leading resonances are larger than they would have been if the coefficient
had been fixed at that value since the beginning.
From a quantitative point of view, the scenario depends both on
the final value and the growth rate of the damping coefficient.
The relevance of the results for the spin-orbit model are discussed in some detail.
\end{abstract}


\zerarcounters
\section{Introduction}
\label{sec:1}

Take a one-dimensional system driven by an external force,
\begin{equation} \label{eq:1.1}
\ddot x + F(x, t ) = 0 
\end{equation}
with $F$ a smooth function $2\pi$-periodic in time $t$;
here and henceforth the dot denotes derivative with respect to time.
Then add a friction term. Usually friction is modelled
as a term proportional to the velocity, so that the equations of motion become
\begin{equation} \label{eq:1.2}
\ddot x + F(x,t) + \g \, \dot x  = 0 , \qquad \g>0 , 
\end{equation}
with the proportionality constant $\g$ referred to as the \emph{damping coefficient}.

As a consequence of friction attractors appear \cite{LL}.
If the system is a perturbation of an integrable one,
there is strong evidence that all attractors are either equilibrium
points (if any) or periodic orbits with periods $T$ which are rational multiples
of the forcing period $2\pi$ \cite{BBDGG}.
If $2\pi/T=p/q$, with $p,q\in\NNN$ and relatively prime, one says that the
periodic orbit is a $p\,$:$\,q$ \emph{resonance}.
For each attractor one can study the corresponding basin of attraction,
that is the set of initial data which approach the attractor as time goes to infinity.
If all motions are bounded, one expects the union of all basins of attraction
to fill the entire phase space, up to a set of zero measure. This appears
to be confirmed by numerical simulations~\cite{FG,RMG,BGG,BBDGG}.

Recently such a scenario has been numerically investigated in several models
of physical interest, such as the dissipative standard map \cite{FG,RMG,FGHY},
the pendulum with oscillating suspension point \cite{BGG},
the driven cubic oscillator \cite{BBDGG}  and the spin-orbit model \cite{Aleshkina,CC}.
What emerges from the numerical simulations is that, for fixed damping,
only a finite number of either point or periodic attractors is present and every
initial datum in phase space is attracted by one of them, according
to Palis' conjecture \cite{P1,P2}. However, which attractors are really
present and the sizes of their basins depend on the value of the damping coefficient. 
If the latter is very small then many periodic attractors can coexist.
This phenomenon is usually called \emph{multistability}; see \cite{FG,FGHY}.
By taking larger values for the damping coefficient, many of the attracting periodic orbits disappear
and, eventually, when the coefficient becomes very large, only a few, if
any, still persist: for every other resonance there is a threshold value for the
damping coefficient above which the corresponding attractor disappears.

In this paper we aim to study what happens when the friction is not fixed but grows in time,
more precisely when the damping coefficient is not a constant but a slowly increasing function of time.
This is a very natural scenario: it is reasonable to suppose that in many physical contexts
dissipation tends to increase to some asymptotic value.
We aim to show (numerically) in such a setting,
with the damping coefficient slowly increasing to a final value,
that the relative areas of some basins of attraction become larger than
they would be if the damping were fixed for all time at the final value.
In other words, we claim that, in order to understand the dynamics of a
forced system in the presence of damping, not only is the final value of the friction important,
but also the time evolution of the damping itself plays a role.
So, by looking in the present at a damped system
which evolved from an original nearly conservative one, with the
friction slowly increasing from virtually zero to the present value,
it can happen that an attractor, which should have a small basin of attraction
on the basis of the final value of the friction, is instead much larger than expected.
Of course, as we shall see, several elements come into the picture, in particular the 
growth rate of the damping coefficient and the closeness between its threshold and final values.

We shall investigate in detail a model system, the driven cubic oscillator
in the presence of friction, which is particularly suited for numerical investigations
because of its simplicity. We shall first study in Section \ref{sec:2} some properties
in the case of constant friction, with some details worked out in Appendix \ref{app:a}.
Then in Section \ref{sec:3} we shall see how the behaviour of the system is affected
by the presence of a non-constant, in fact slowly increasing friction. In Section \ref{sec:4}
we shall introduce another system of physical interest, the spin-orbit system,
with some details deferred to Appendices \ref{app:b} to \ref{app:e},
and we shall discuss how the results described in the previous sections
may be relevant to the study of its behaviour. Further comments are deferred to
Section \ref{sec:5}. Finally in Section \ref{sec:6} we draw our
conclusions and briefly discuss open problems. Some discussion
of the codes we used for the numerical analysis is given in Appendix \ref{app:f}. 

\zerarcounters
\section{The driven cubic oscillator with constant friction}
\label{sec:2}

Let us consider the cubic oscillator, subject to periodic forcing and
in the presence of friction,
\begin{equation} \label{eq:2.1}
\ddot x + x^{3} + \e \, f(t) \, x^{3} + \g \, \dot x = 0 , \qquad f(t)=\cos t ,
\end{equation}
where $x\in\RRR$ and $\e$ is a real parameter, called the \emph{perturbation parameter},
that we shall suppose positive (for definiteness). Of course one could
consider more general expressions for $f$ and the choice made here is for
simplicity.  The system \eqref{eq:2.1}
has been investigated in \cite{BBDGG}, with $\g$ a fixed positive constant.
The constants $\e$ and $\g$ are two control parameters, measuring
respectively the forcing and the dissipation of the system.

We shall look at \eqref{eq:2.1} as a non-autonomous first order differential equation,
so that the phase space is $\RRR^{2}$. Note that $(x,\dot x)=(0,0)$ is
an equilibrium point for all values of $\e$ and $\g$. Moreover, for $\e=\g=0$
the system is integrable and all motions are periodic. One can write
the solutions explicitly in terms of elliptic integrals \cite{GR,BBDGG}.
For $\e\neq 0$ (hence $\e>0$) and $\g>0$ fixed,
a finite number of periodic orbits of the unperturbed system persist and, together with the
equilibrium point, they attract every trajectory in phase space \cite{BBDGG}.
Such periodic orbits are called \emph{subharmonic solutions} in the literature \cite{GH}.
Each periodic orbit can be identified through the corresponding frequency or, better, the
ratio $\om:=p/q$ between its frequency and the frequency of the forcing term.
For each periodic orbit one can compute the corresponding \emph{threshold value}
$\g(\om,\e)$: if $\g>\g(\om,\e)$ the orbit ceases to exist, while
for $\g<\g(\om,\e)$ the orbit is present, with a basin of attraction whose area
depends on the actual value of $\g$.  At fixed $\e$, one has
\begin{equation} \label{eq:2.2}
\lim_{\max\{p,q\}\to \infty} \g(p/q,\e) = 0 .
\end{equation}
Therefore, if we assume that all attractors different from the equilibrium point are periodic
and no periodic attractors other than subharmonic solutions exist,
then we find that at fixed $\e$ and $\g$ only a finite number of attractors exists.
We note that the assumption above, even though we have no proof,
is consistent with numerical findings \cite{BBDGG}.

The threshold value $\g(\om,\e)$ depends smoothly on $\e$ \cite{CH,GH,GBD}:
for all $\om\in\QQQ$ there exists $n(\om)\in\NNN$ such that the corresponding threshold value
is of the form $\g(\om,\e)=C_{0}(\om,\e)\,\e^{n(\om)}$,  with the constant $C_{0}(\om,\e)$
nearly independent of $\e$ for $\e$ small; more precisely $C_{0}(\om,\e)$
tends to a constant $C_{0}(\om)$ as $\e$ goes to zero, so that we can consider it a constant
for $\e$ small enough. Resonances are classified as follows:
we refer to resonances with frequency $\om$ such that $n(\om)=1$ as \emph{primary},
to resonances with frequency $\om$ such that $n(\om)=2$ as \emph{secondary}, 
and so on \cite{PC,BBDGG}. Of course such a classification makes sense only for $\e$ small enough.
The primary resonances are the most important, in the sense that,
at fixed small $\e$, for $\g$ large enough, only primary resonances are present;
moreover, by decreasing the value of $\g$, although non-primary resonances
appear, they have a small basin of attraction with respect to those of the primary ones.
The threshold values of the leading attractors (that is the attractors
with largest threshold values), in terms of the constants $C_{0}(\om)$,
were computed analytically in \cite{BBDGG} and are reproduced in
Tables \ref{tab:2.1} and \ref{tab:2.2}.  In particular the periodic attractors
with frequency $1/q$, with $q$ odd, appear in pairs \cite{BBDGG}.
The higher order corrections to $C_{0}(\om,\e)$ are explicitly computable; 
however we shall not need to do this here. Note that the classification of resonances
and the corresponding threshold values strongly depend on the forcing:
all values in this and next Section refer to $f(t)=\cos t$, as in (\ref{eq:2.1}).

\begin{table}[ht]
\centering
\caption{Values of the constants $C_{0}(p/q)$ for $p=1$ and $q=2,4,6,8,10$
(leading primary resonances) for the cubic oscillator \eqref{eq:2.1}; the threshold
values are of the form $\g(\om,\e)=C_{0}(\om)\e+O(\e^{2})$.}
\begin{center}
\setlength\tabcolsep{5pt}
\vskip-.5truecm
\vrule
\begin{tabular}{llllllllllll}
\hline\noalign{\smallskip}\hline
$q$ & \vrule\vrule\vrule & $2$ & \vrule & $4$ & \vrule & $6$ & \vrule & $8$ & \vrule & $10$ \\
\hline\noalign{\smallskip}\hline
$C_{0}(1/q)$ & \vrule\vrule\vrule & $0.178442$ & \vrule & $0.061574$ & \vrule & $0.008980$
& \vrule & $0.000920$ & \vrule & $0.000078$ \\
\hline
\end{tabular}
\hspace{-0.1cm}\vrule
\vskip-.5truecm
\end{center}
\label{tab:2.1}
\end{table}

\begin{table}[ht]
\centering
\caption{Values of the constants $C_{0}(p/q)$ for $p=1$ and $q=1,3,5,7,9$
(leading secondary resonances) for the cubic oscillator \eqref{eq:2.1}; the threshold
values are of the form $\g(\om,\e)=C_{0}(\om)\e^{2}+O(\e^{3})$.}
\begin{center}
\setlength\tabcolsep{5pt}
\vskip-.5truecm
\vrule
\begin{tabular}{llllllllllll}
\hline\noalign{\smallskip}\hline
$q$ & \vrule\vrule\vrule & $1$ & \vrule & $3$ & \vrule & $5$ & \vrule & $7$ & \vrule & $9$ \\
\hline\noalign{\smallskip}\hline
$C_{0}(1/q)$ & \vrule\vrule\vrule & $0.146322$ & \vrule & $0.065001$ & \vrule & $0.006488$
& \vrule & $0.000177$ & \vrule & $0.000002$ \\
\hline
\end{tabular}
\hspace{-0.1cm}\vrule
\vskip-.5truecm
\end{center}
\label{tab:2.2}
\end{table}

Consider the system \eqref{eq:2.1} at fixed $\e$. 
For $\g$ large enough, the only attractor left is the equilibrium point;
in that case all trajectories eventually go toward this point, which
becomes a global attractor (see Appendix \ref{app:a}).
If $\g$ is not too large --- that is, according to Table \ref{tab:2.1},
if $\g<C_{0}(1/2)\e$, up to higher order corrections ---
then, besides the equilibrium point, there is a finite number of other attractors,
which are periodic orbits.  

For $\e=0.1$, from Tables \ref{tab:2.1} and \ref{tab:2.2} one obtains the threshold values
$\g(1/2,0.1)\approx 0.018$, $\g(1/4,0.1) \approx 0.0062$, $\g(1,0.1) \approx 0.0015$,
$\g(1/6,0.1) \approx 0.00090$, $\g(1/3,0.1) \approx 0.00065$, $\g(1/8,0.1) \approx 0.000092$,
$\g(1/5,0.1) \approx 0.000065$, and so on.
Take the initial data in a finite domain of phase space, say the square $\calQ=[-1,1]\times[-1,1]$:
the relative areas of the parts of the basins of attraction contained
in $\calQ$ for some values of $\g$ are given in Table \ref{tab:2.3}.
In principle the relative areas depend on the domain, but one expects
that they do not change too much by changing the domain, provided the latter is not too small.
Note that for $\e=0.1$ and $\g \le 0.00005$, other attractors than those listed
in Table \ref{tab:2.3} appear (namely periodic orbits with frequencies
$1/8$, $1/5$ and $3/4$ for $\g=0.00005$;
with frequencies
$1/8$, $1/5$, $3/10$, $2/5$, $5/12$ and $3/4$ for $\g=0.00001$;
and with frequencies
$1/10$, $1/8$, $1/7$, $1/5$, $3/14$, $2/7$, $3/10$, $2/5$, $5/12$,
$3/7$, $2/3$ and $3/4$ for $\g=0.000005$),
so explaining why the relative areas of the
basins of attractions considered there do not sum up to 100\%.
Small discrepancies for the other values are simply due to round-off error
(the error on the data is in the first decimal digit; see Appendix \ref{app:f}).

\begin{table}[ht]
\centering
\caption{Relative areas $A(\omega, \gamma)$, \%, of the parts of the basins of attraction contained inside
the square $\calQ$ for $\e=0.1$ and some values of $\g$. The attractors are identified
by the corresponding frequency (0 is the origin). 
The number of random initial conditions taken in $\calQ$ is
$1\,000\,000$ up to $\g=0.0001$, $500\,000$ for $\g=0.00005$,
$150\,000$ for $\g=0.00001$ and $50\,000$ for $\g=0.000005$.}
\begin{center}
\setlength\tabcolsep{5pt}
\vskip-.5truecm
\vrule
\begin{tabular}{llllllllllllllllllllll}
\hline\noalign{\smallskip}\hline
$\om$ & \vrule\vrule\vrule & 0  & \vrule &
1/2 & \vrule & 1/4 & \vrule & 1a  & \vrule & 1b 
& \vrule & 1/6  & \vrule & 1/3a  & \vrule & 1/3b & \vrule & 3/8 \\        
\hline\noalign{\smallskip}\hline
$\g=0.020000$ & \vrule\vrule\vrule & $\!\!\!100.0$ & \vrule
& $00.0$ & \vrule & $00.0$ & \vrule & $00.0$ & \vrule
& $00.0$ & \vrule & $00.0$ & \vrule & $00.0$ & \vrule & $00.0$ & \vrule & $00.0$\\
\hline
$\g=0.015000$ & \vrule\vrule\vrule & $91.1$ & \vrule
& $08.9$ & \vrule & $00.0$ & \vrule & $00.0$ & \vrule
& $00.0$ & \vrule & $00.0$ & \vrule & $00.0$ & \vrule & $00.0$ & \vrule & $00.0$ \\
\hline
$\g=0.010000$ & \vrule\vrule\vrule & $79.1$ & \vrule
& $20.9$ & \vrule & $00.0$ & \vrule & $00.0$ & \vrule
& $00.0$ & \vrule & $00.0$ & \vrule & $00.0$ & \vrule & $00.0$ & \vrule & $00.0$ \\
\hline
$\g=0.005000$ & \vrule\vrule\vrule & $64.9$ & \vrule
& $31.8$ & \vrule & $03.4$ & \vrule & $00.0$ & \vrule
& $00.0$ & \vrule & $00.0$ & \vrule & $00.0$ & \vrule & $00.0$ & \vrule & $00.0$ \\
\hline
$\g=0.001000 $& \vrule\vrule\vrule & $44.5$ & \vrule
& $40.9$ & \vrule & $13.2$ & \vrule & $00.7$ & \vrule
& $00.7$ & \vrule & $00.0$ & \vrule & $00.0$ & \vrule & $00.0$ & \vrule & $00.0$ \\
\hline
$\g=0.000500$ & \vrule\vrule\vrule & $38.7$ & \vrule
& $41.8$ & \vrule & $14.7$ & \vrule & $01.3$ & \vrule
& $01.3$ & \vrule & $01.7$ & \vrule & $00.3$ & \vrule & $00.3$ & \vrule & $00.0$ \\
\hline
$\g=0.000100$ & \vrule\vrule\vrule & $32.2$ & \vrule
& $41.9$ & \vrule & $14.0$ & \vrule & $02.6$ & \vrule
& $02.6$ & \vrule & $03.6$ & \vrule & $01.5$ & \vrule & $01.5$ & \vrule & $00.1$ \\
\hline
$\g=0.000050$ & \vrule\vrule\vrule & $30.2$ & \vrule
& $41.6$ & \vrule & $13.8$ & \vrule & $02.8$ & \vrule
& $02.8$ & \vrule & $03.8$ & \vrule & $01.7$ & \vrule & $01.7$ & \vrule & $00.6$ \\
\hline
$\g=0.000010$ & \vrule\vrule\vrule & $26.9$ & \vrule
& $41.1$ & \vrule & $13.2$ & \vrule & $02.9$ & \vrule
& $02.9$ & \vrule & $03.9$ & \vrule & $01.8$ & \vrule & $01.8$ & \vrule & $01.1$ \\
\hline
$\g=0.000005$ & \vrule\vrule\vrule & $26.2$ & \vrule
& $40.9$ & \vrule & $13.0$ & \vrule & $02.9$ & \vrule
& $02.9$ & \vrule & $03.8$ & \vrule & $01.8$ & \vrule & $01.8$ & \vrule & $01.3$ \\\hline
\end{tabular}
\hspace{-0.1cm}\vrule
\vskip-.5truecm
\end{center}
\label{tab:2.3}
\end{table}

\begin{figure}[ht] 
\centering 
\ifthenelse{\boolean{pdflat}}
{\subfigure{\includegraphics*[width=2.8in]{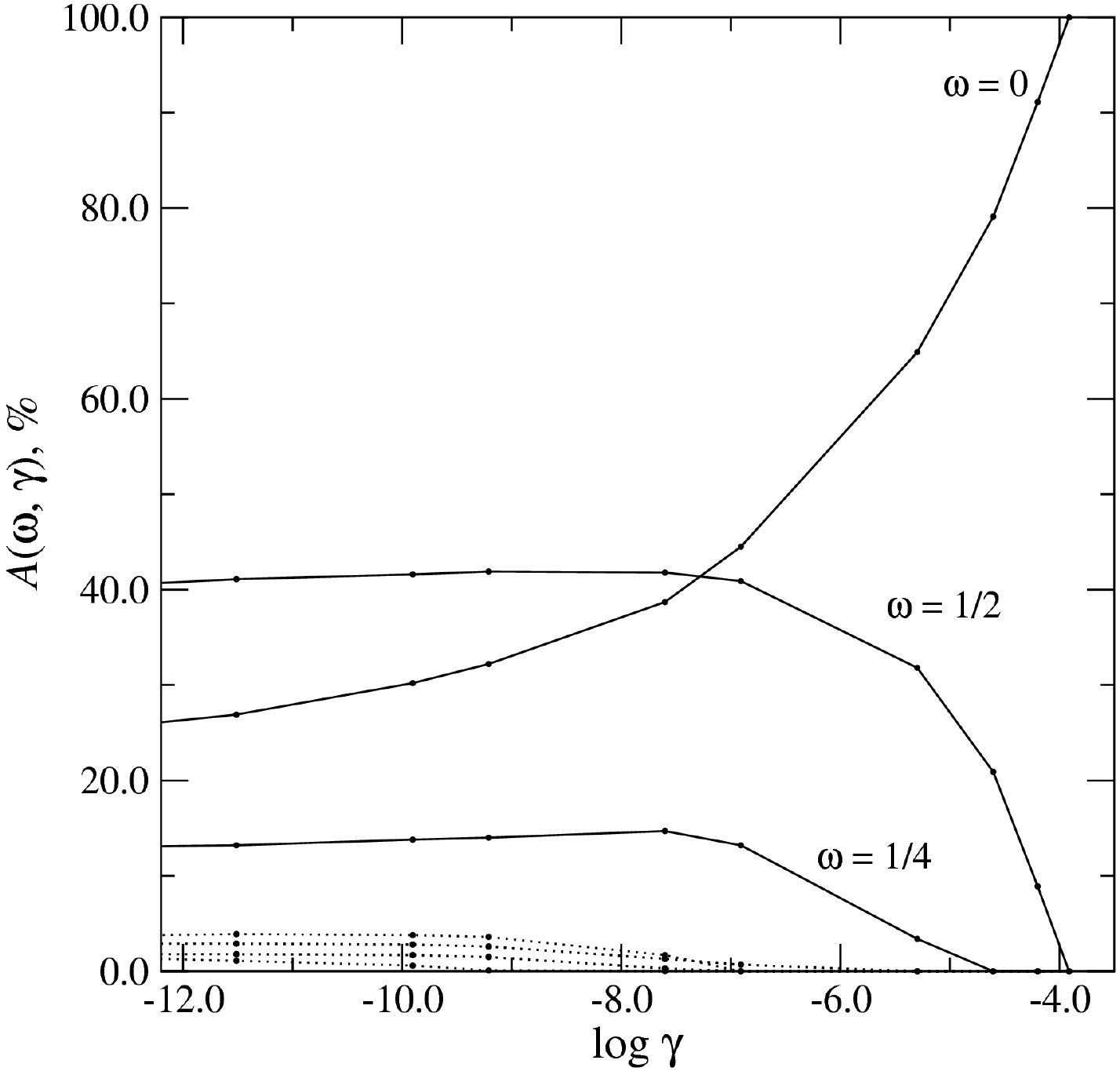}
\hskip2.truecm
\subfigure{\includegraphics*[width=2.8in]{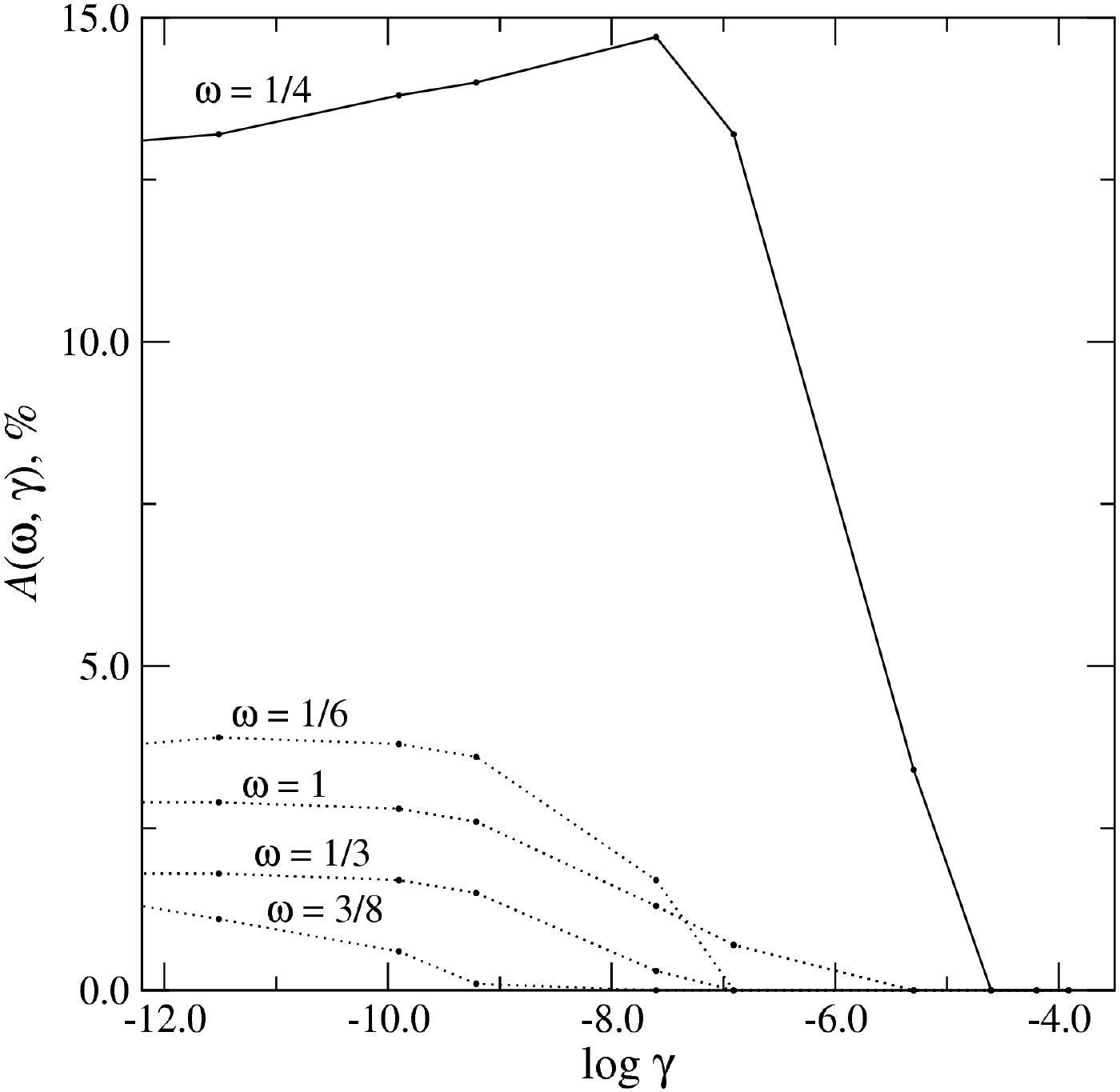} }}}
{\subfigure{\includegraphics*[width=2.8in]{fig21.eps}
\hskip2.truecm
\subfigure{\includegraphics*[width=2.8in]{fig22.eps} }}}
\caption{Relative areas $A(\om,\g)$ of the basins of attraction versus $\log \g$
for the values of $\g$ listed in Table \ref{tab:2.3} (right-hand figure) and
a magnification for the periodic orbits
with $\om=1/4,1,1/6,1/3,3/8$ (left-hand figure).}
\label{fig:2.1} 
\end{figure}

If one plots the relative areas $A(\om,\g)$ of the basins of attraction 
versus $\g$ one finds the situation depicted in Figure \ref{fig:2.1}.
Of course in general $A(\om,\g)$ depends also on $\e$, i.e. $A(\om,\g)=A(\om,\g,\e)$,
although we are not making explicit such a dependence since $\e$ has been fixed
at $\e=0.1$; same comment applies to the quantities $A_{\rm max}(\om)=A_{\rm max}(\om,\e)$
and $A(\om)=A(\om,\e)$ to be introduced.

It can be seen that for any $\om\in\QQQ$ one has $A(\om,\g)=0$ if $\g>\g(\om,\e)$.
By decreasing $\g$ below $\g(\om,\e)$, $A(\om,\g)$ increases up to a maximum value
$A_{\rm max}(\om)$ which tends to stabilise. For very small values of $\g$
one observes a slight bending downward. It would be interesting to
investigate very --- in principle arbitrarily --- small values of $\g$,
but of course we have to cope with the technical limitations of computation:
studying arbitrarily small friction would require running programs
for arbitrarily large times and with arbitrarily high precision --- see
also comments in Appendix \ref{app:f}.
However, by looking at Figure \ref{fig:2.1} and noting that numerical
evidence suggests that all attractors different from the origin are periodic orbits,
we make the following conjecture: as $\g$ goes to $0^{+}$, for all $\om\in\QQQ$
the relative area $A(\om,\g)$ tends to a finite limit $A(\om)$ such that
\begin{equation} \label{eq:2.3}
\sum_{\om\in\QQQ} A(\om) = 100\% ,
\end{equation}
where $\om=0$ designates the origin. Of course when $\g=0$ the area of each basin drops to zero, 
so that, accepting the conjecture above, all functions $A(\om,\g)$ are discontinuous at $\g=0$.
This is not surprising: a similar situation arises in the absence of forcing,
where the only attractor is the origin, with a basin of attraction
which passes abruptly from zero ($\g=0$) to 100\% ($\g>0$). Moreover,
analogously to \cite{RMG}, we would expect that the total number of
periodic attractors $N_{\rm p}$ grows as a power of $\g$ when $\g$ tends to $0$.
This means that if the limits $A(\om)$ vanished at least one function $A(\om,\g)$
should be exponential in $\log\g$, a behaviour which seems
unlikely in the light of Figure \ref{fig:2.1}.



\zerarcounters
\section{The driven cubic oscillator with increasing friction}
\label{sec:3}

Here we shall consider $\g=\g(t)$ explicitly depending on time, that is
\begin{equation} \label{eq:3.1}
\ddot x + x^{3} + \e \, f(t) \, x^{3} + \g(t) \, \dot x = 0 , \qquad f(t)=\cos t ,
\end{equation}
For both concreteness and simplicity reasons, we shall consider a
dissipation $\g(t)$ linearly increasing in time up to some final value, i.e
\begin{equation} \label{eq:3.2}
\gamma(t) = \begin{cases}
\displaystyle{\gamma_{0} \frac{t}{T_{0}}} , & \qquad 0 \le t < T_{0} , \\
\gamma_{0} , & \qquad t \ge T_{0} , \end{cases}
\end{equation}
where the parameters $\gamma_{0}$ and $T_{0}$ are positive constants.
However, the results we are going to describe should not depend too much
on the exact form of the function $\g(t)$, as long as it is a slowly
increasing function; see Section \ref{sec:5}.
In (\ref{eq:3.2}) we shall take $T_{0}=\Delta/\g_{0}$, whose form is
suggested by the fact that trajectories converge toward
an attractor at a rate proportional to $1/\g_{0}$ (see Appendix \ref{app:a}).

\begin{figure}[ht] 
\centering 
\ifthenelse{\boolean{pdflat}}
{\includegraphics*[width=2.8in]{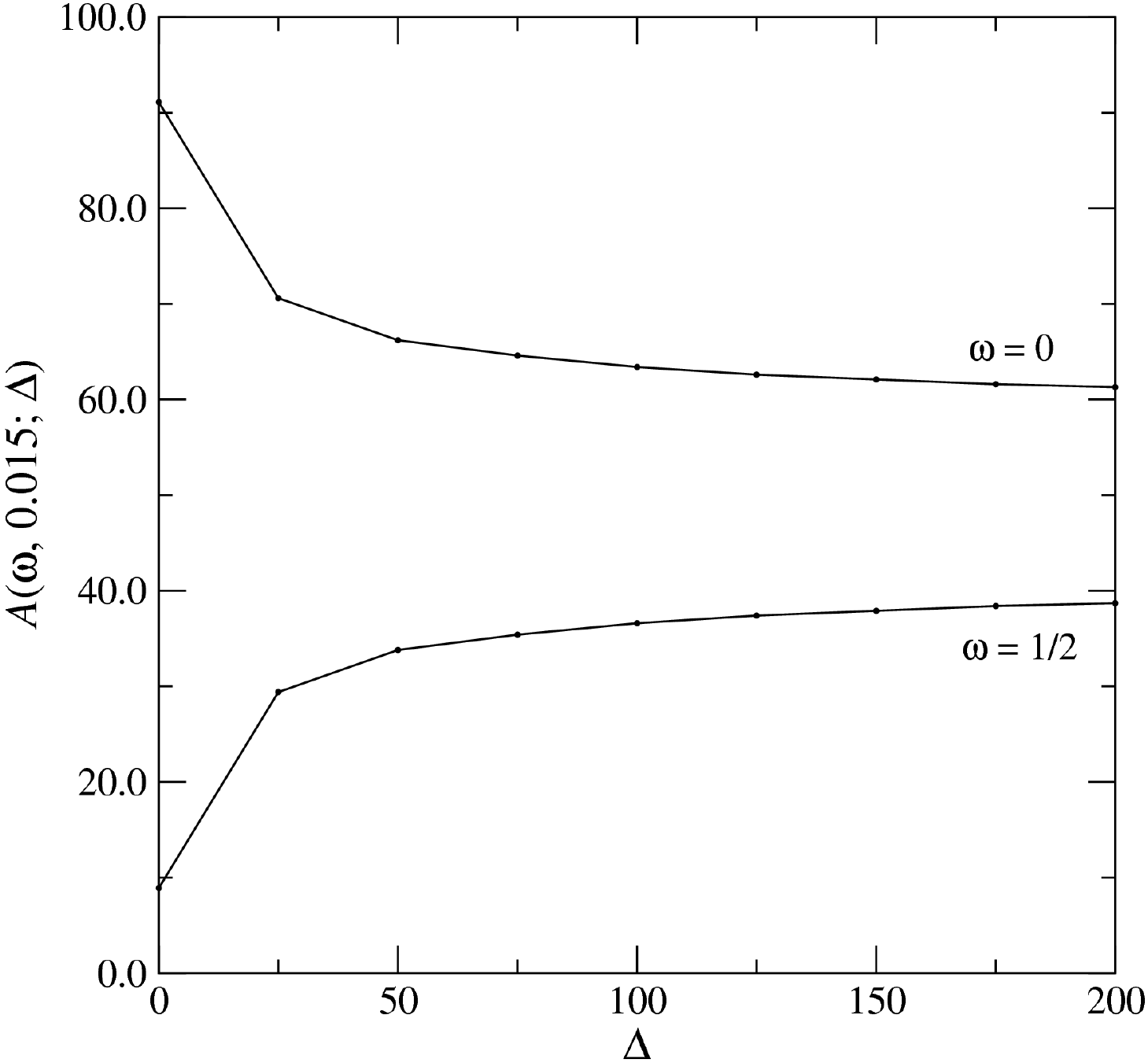} }
{\includegraphics*[width=2.8in]{fig31.eps} }
\caption{Relative measures of the basins of attraction versus $\Delta$ for $\g_{0}=0.015$.}
\label{fig:3.1} 
\end{figure}

Hence, consider the system with $\e=0.1$ again but now with $\g(t)$
given by \eqref{eq:3.2}, with $T_{0}=\Delta/\g_{0}$ and $\g_{0}=0.015$.
Computing the corresponding relative areas $A(\om,\g_{0};\Delta)=A(\om,\g_{0},0.1;\Delta)$
--- that is $A(\om,\g_{0},\e;\Delta)$ for $\e=0.1$ ---
for different values of $\Delta$, we obtain the results in Table \ref{tab:3.1} and Figure \ref{fig:3.1}.
If $\Delta$ is very small, the damping coefficient reaches the asymptotic
value $\g_0$ almost immediately, and we would expect to obtain the same scenario
as in the previous case ($\g$ constant): two attractors, corresponding to the origin
and the 1:2 resonance,  with basins whose relative areas are close to
the values for $\Delta=0$, i.e. 91.1\% and 8.9\%, respectively. 
On the other hand, if $\Delta$ becomes larger, we find that
the relative area of the basin of attraction of the origin decreases, whereas that
of the basin of attraction of the 1:2 resonance increases.
For $\Delta$ very large, these areas apparently tend to constant values
of around $61\%$ and $39\%$, respectively; see Table \ref{tab:3.1}. 

\begin{table}[ht]
\centering
\caption{Relative areas $A(\om,0.015;\Delta)$ of the parts of the basins of attraction
contained inside the square $\calQ$
for $\e=0.1$ and $\g(t)$ given by (\ref{eq:3.2}) with $\g_{0}=0.015$
and $T_{0}=\Delta/\g_{0}$, for various values of $\Delta$ and $\om=0,1/2$
($\om=0$ is the origin). In each case, $1\,000\,000$ random initial conditions have been taken in $\calQ$.}
\begin{center}
\setlength\tabcolsep{5pt}
\vskip-.5truecm
\vrule
\begin{tabular}{lllllllllllllllllllll}
\hline\noalign{\smallskip}\hline
$\Delta$ & \vrule\vrule\vrule & $0$ & \vrule & $25$ & \vrule & $50$
& \vrule & $75$ & \vrule & $100$ & \vrule & $125$ & \vrule & $150$
& \vrule & $175$ & \vrule & $200$ \\
\hline\noalign{\smallskip}\hline
   $\om=0$ & \vrule\vrule\vrule & $91.1$ & \vrule & $70.6$ & \vrule & $66.2$ & \vrule & $64.6$ 
& \vrule & $63.4$ & \vrule & $62.6$ & \vrule & $62.1$ &\vrule & $61.6$ & \vrule & $61.3$ \\
\hline
$\om=1/2$ & \vrule\vrule\vrule & $08.9$ & \vrule & $29.4$ & \vrule & $33.8$ & \vrule & $35.4$ 
& \vrule & $36.6$ & \vrule & $37.4$ & \vrule & $37.9$ & \vrule & $38.4$ & \vrule & $38.7$ & \\
\hline
\end{tabular}
\hspace{-0.1cm}\vrule
\vskip-.5truecm
\end{center}
\label{tab:3.1}
\end{table}

\begin{figure}[ht] 
\centering 
\ifthenelse{\boolean{pdflat}}
{\includegraphics*[width=2.8in]{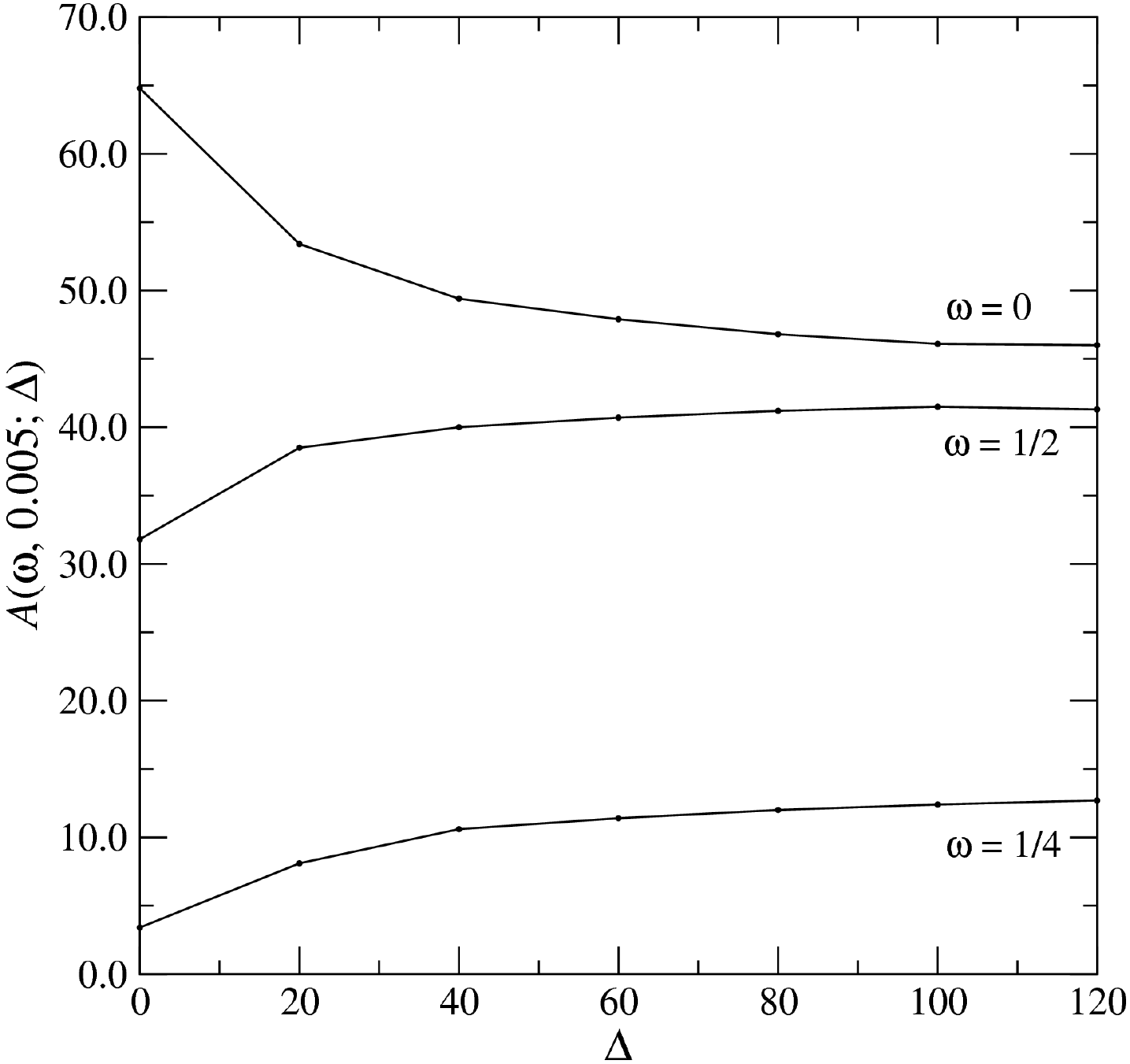} }
{\includegraphics*[width=2.8in]{fig32.eps} }
\caption{Relative areas of the basins of attraction versus $\Delta$ for $\g_{0}=0.005$.}
\label{fig:3.2} 
\end{figure}

Analogous results are found, for instance, for $\g_{0}=0.005$;
see Table \ref{tab:3.2} and Figure \ref{fig:3.2}. One sees that for $\Delta=20$
the relative areas of the basins of attraction of the origin and of the
1:2 and 1:4 resonances have already appreciably changed: they have become,
respectively, 53.4\%, 38.5\% and 8.1\%. By further increasing $\Delta$, once again
a saturation phenomenon is observed and the relative areas settle
about asymptotic values around 45\%, 42\% and 13\% (for instance for $\Delta=120$
the areas are, respectively, 46.0\%, 41.3\% and 12.7\%).
Note that the value $\g_{0}=0.005$ is such that the threshold values
$\g(1/2,0.1) \approx 0.018$ and $\g(1/4,0.1)\approx 0.0062$ of the persisting
resonances are slightly above it (that is their ratios with $\g_{0}$ are of order 1).

\begin{table}[ht]
\centering
\caption{Relative areas $A(\om,0.005;\Delta)$ 
of the parts of the basins of attraction contained inside the square $\calQ$
for $\e=0.1$ and $\g(t)$ given by (\ref{eq:3.2}) with $\g_{0}=0.005$
and $T_{0}=\Delta/\g_{0}$, for various values of $\Delta$ and $\om=0,1/2,1/4$
($\om=0$ is the origin). $500\,000$ random initial conditions have been taken in $\calQ$.}
\begin{center}
\setlength\tabcolsep{5pt}
\vskip-.5truecm
\vrule
\begin{tabular}{lllllllllllllllll}
\hline\noalign{\smallskip}\hline
$\Delta$ & \vrule\vrule\vrule & $0$ & \vrule & $20$ & \vrule & $40$
& \vrule & $60$ & \vrule & $80$ & \vrule & $100$ & \vrule & $120$ \\
\hline\noalign{\smallskip}\hline
   $\om = 0$ & \vrule\vrule\vrule & $64.8$ & \vrule & $53.4$ & \vrule & $49.4$ 
&\vrule & $47.9$ & \vrule & $46.8$ & \vrule & $46.1$ & \vrule & $46.0$ \\
\hline
$\om = 1/2$ & \vrule\vrule\vrule & $31.8$ & \vrule & $38.5$ & \vrule & $40.0$ 
&\vrule & $40.7$ & \vrule & $41.2$ & \vrule & $41.5$ & \vrule & $41.3$ \\ 
\hline
$\om = 1/4$ & \vrule\vrule\vrule & $03.4$ & \vrule & $08.1$ & \vrule & $10.6$
&\vrule & $11.4$ & \vrule & $12.0$ & \vrule & $12.4$ & \vrule & $12.7$ \\ 
\hline
\end{tabular}
\hspace{-0.1cm}\vrule
\vskip-.5truecm
\end{center}
\label{tab:3.2}
\end{table}

\begin{figure}[ht] 
\centering 
\ifthenelse{\boolean{pdflat}}
{\includegraphics*[width=2.8in]{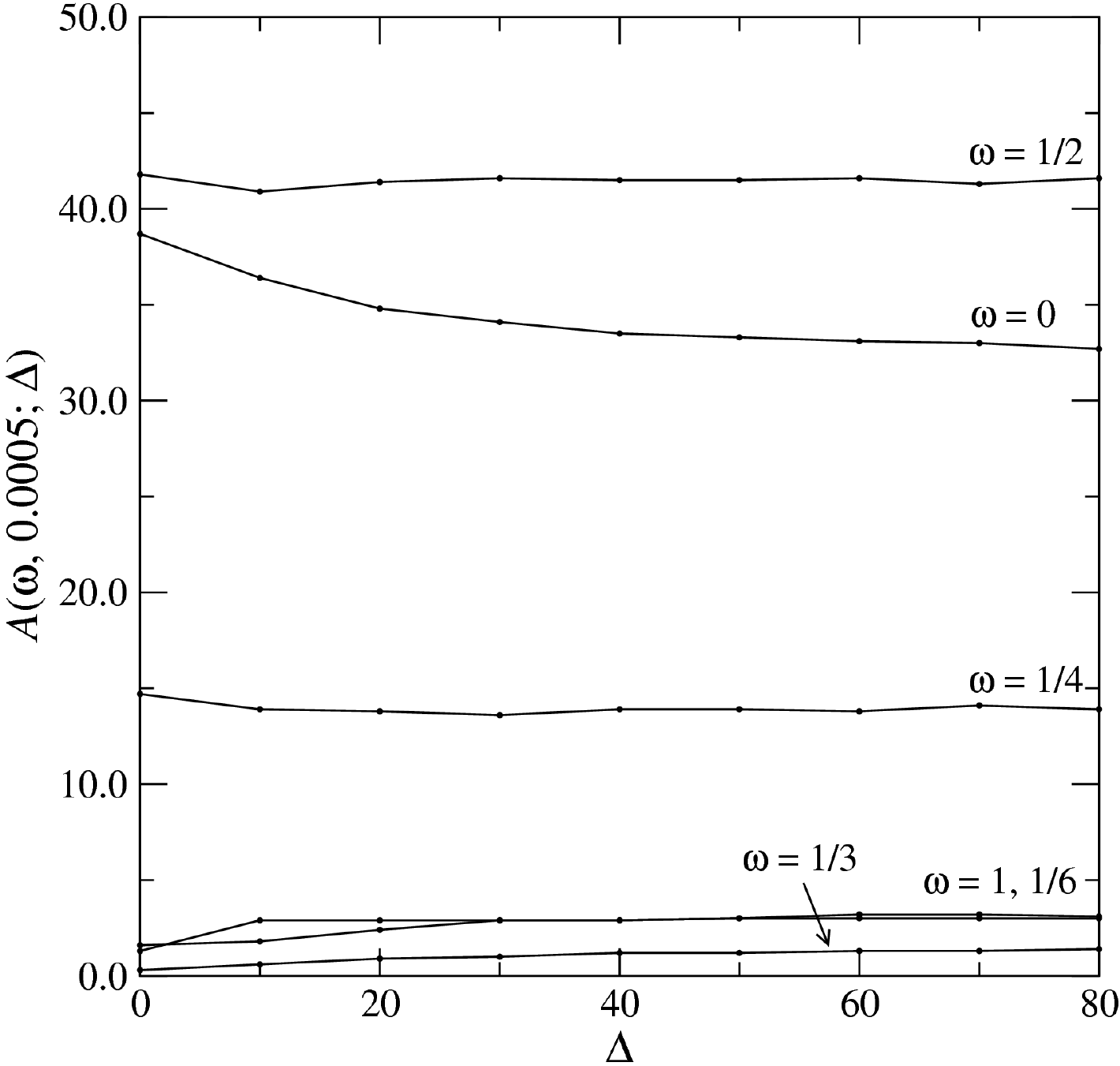} }
{\includegraphics*[width=2.8in]{fig33.eps} }
\caption{Relative areas of the basins of attraction versus $\Delta$ for $\g_{0}=0.0005$.}
\label{fig:3.3} 
\end{figure}

If one fixes the value $\g_{0}=0.0005$,
then the 1:6, 1:1 and 1:3 resonances are also present. On the other hand
the threshold values  of the 1:2 and 1:4 resonances are appreciably larger then $\g_{0}$
(that is $\g(\om) \gg \g_{0}$ for $\om=1/2$ and $\om=1/4$),
whereas the threshold values $\g(1/6,0.1) \approx 0.00090$,
$\g(1,0.1)\approx 0.0015$ and $\g(1/3,0.1) \approx 0.00065$
of the 1:6, 1:1 and 1:3 resonances, respectively, are not too different from $\g_{0}$.
If we again take $\g(t)$ as in \eqref{eq:3.2}, with $T_{0}=\Delta/\g_{0}$
and $\g_{0}=0.0005$, we have the results in Table \ref{tab:3.3} and Figure \ref{fig:3.3}.

\begin{table}[ht]
\centering
\caption{Relative areas $A(\om,0.0005;\Delta)$  
of the parts of the basins of attraction contained inside $\calQ$
for $\e=0.1$ and $\g(t)$ given by (\ref{eq:3.2}) with $\g_{0}=0.0005$
and $T_{0}=\Delta/\g_{0}$, for various values of $\Delta$
and $\om=0,1/2,1/4,1,1/6,1/3$ ($\om=0$ is the origin).
$250\,000$ random initial conditions have been taken in $\calQ$.}
\begin{center}
\setlength\tabcolsep{5pt}
\vskip-.5truecm
\vrule
\begin{tabular}{llllllllllllllllllll}
\hline\noalign{\smallskip}\hline
$\Delta$ & \vrule\vrule\vrule & $0$ & \vrule & $10$ & \vrule & $20$
& \vrule & $30$ & \vrule & $40$ & \vrule & $50$ 
& \vrule & $60$ & \vrule & $70$ & \vrule & $80$ \\
\hline\noalign{\smallskip}\hline
  $\om = 0$ & \vrule\vrule\vrule & $38.7$ & \vrule & $36.4$ & \vrule & $34.8$ 
&\vrule & $34.1$ & \vrule & $33.5$ & \vrule & $33.3$ 
&\vrule & $33.1$ & \vrule & $32.9$ & \vrule & $32.7$ \\ 
\hline
$\om = 1/2$ & \vrule\vrule\vrule & $41.8$ & \vrule & $40.9$ & \vrule & $41.4$ 
&\vrule & $41.6$ & \vrule & $41.5$ & \vrule & $41.5$ 
&\vrule & $41.6$ & \vrule & $41.3$ & \vrule & $41.6$ \\
\hline 
$\om = 1/4$ & \vrule\vrule\vrule & $14.7$ & \vrule & $13.9$ & \vrule & $13.8$
&\vrule & $13.6$ & \vrule & $13.9$ & \vrule & $13.9$ 
&\vrule & $13.8$ & \vrule & $14.0$ & \vrule & $13.9$ \\
\hline
$\om = 1$ & \vrule\vrule\vrule & $01.3$ & \vrule & $02.9$ & \vrule & $02.9$ 
&\vrule & $02.9$ & \vrule & $02.9$ & \vrule & $03.0$ 
&\vrule & $03.0$ & \vrule & $03.0$ & \vrule & $03.0$ \\
\hline
$\om = 1/6$ & \vrule\vrule\vrule & $01.6$ & \vrule & $01.8$ & \vrule & $02.4$ 
&\vrule & $02.9$ & \vrule & $02.9$ & \vrule & $03.0$ 
&\vrule & $03.2$ & \vrule & $03.1$ & \vrule & $03.2$ \\
\hline
$\om = 1/3$ & \vrule\vrule\vrule & $00.3$ & \vrule & $00.6$ & \vrule & $00.9$
&\vrule & $01.0$ & \vrule & $01.2$ & \vrule & $01.2$ 
&\vrule & $01.3$ & \vrule & $01.4$ & \vrule & $01.3$ \\
\hline
\end{tabular}
\hspace{-0.1cm}\vrule
\vskip-.5truecm
\end{center}
\label{tab:3.3}
\end{table}

Therefore we have, from a qualitative point of view, the same scenario as in the
case $\g_{0}=0.005$, but with some relevant quantitative differences:
the relative areas of the basins of the 1:2 and 1:4 resonances 
are not too different in the two situations $\g$ constant and $\g$ increasing.

We now give an argument to explain why the basins of attraction are different if
$\g$ is not fixed ab initio to some value $\g_{0}$ but slowly tends to that value.
According to Table \ref{tab:2.3}, for smaller values of $\g$ the
basins of the periodic  attractors are larger. For instance
for $\g=0.005$ the 1:2 and 1:4 resonances have basins with relative areas
31.8\% and 3.4\%, respectively, while the basins of attraction of the same resonances 
for $\g=0.0005$ have relative areas 41.8\% and 14.7\%, respectively.
Then, if we suppose that the friction is slowly increasing in time,
when it passes, say, from $\g=0.0005$ to $0.005$, on the one hand
the size of the basin would decrease because of the larger value of $\g$,
but on the other hand many trajectories have already nearly reached the basin
and hence continue to be attracted toward that resonance.
If friction increases slowly enough we can assume that it is quasi-static.
Therefore, at every instant $\tau$, the
basin of attraction of any resonance has the size corresponding
to the value $\g(\tau)$ at that instant, as can be deduced 
by interpolation from Table \ref{tab:2.3} (or Figure \ref{fig:2.2}),
while the rate of approach to the resonance
can be roughly estimated as proportional to $1/\g(\tau)$;
see Appendix \ref{app:a}. Therefore if $\Delta$ is large enough
(that is if the growth of $\g(t)$ is slow enough)
one expects the trajectory to be captured by
the resonance faster than how the basin of attraction is decreasing.

By increasing the friction further, the basin of a resonance $\om$
can become negligible, until the resonance itself disappears.
If this does not happen, that is if $\g(\om,0.1)>\g_{0}$, then there
is a value of $\Delta$ above which  the relative measure $A(\om,\g_{0};\Delta)$
of the basin saturates to a value close to the maximum value $A_{\rm max}(\om)$
(possibly a bit smaller because of the slight bending downward observed
in Figure \ref{fig:2.1}).
In particular, this explains the difference between  Figures \ref{fig:3.2} and \ref{fig:3.3}.
For concreteness let us focus on the 1:2 resonance.
With respect to the case $\g_{0}=0.0005$, according to Figure \ref{fig:3.1},
the relative area $A(1/2,\g)$  of the basin of attraction does not increase appreciably
when taking smaller values $\g < \g_{0}$: indeed $A(1/2,0.0005)$ is already
close to $A_{\rm max}(1/2)$.

We conclude that the main effect of friction slowly growing to
a final value $\g_{0}$, is that eventually every basin of attraction has essentially
the same size that would appear for lower values of friction.
So, if the basin of attraction of any $p:q$ resonance is larger 
for values of friction lower than the final value,
then, when the final value $\g_{0}$ is reached, one observes a
basin of attraction with relative area larger than $A(p/q,\g_{0})$.
If on the contrary it is more or less the same, then one observes
essentially the same basin one would have by taking
the friction fixed at that value since the beginning.
In other words, if the friction increases in time, one can really
have a larger basin of attraction only if the final value $\g_{0}$
is close enough to the threshold value $\g(\om,0.1)$. However,
if $\g_{0}$ is too close to $\g(\om,0.1)$, for the phenomenon to really occur, 
the rate of growth has to be slow enough: the closer $\g_{0}$ is to
$\g(\om,0.1)$, the larger $\Delta$ to be chosen. For instance,
for $\e=0.1$ and $\g_{0}=0.0005$,  a glance at Table \ref{tab:3.3}
gives the following picture. The areas of the basins of attraction
for the 1:2 and 1:4 resonances have small variations for different values of $\Delta$,
whereas the areas of the basins of attraction for 1:1, 1:3 and 1:6 change
in a more appreciable way when $\Delta$ becomes larger.
Moreover, the threshold value $\g(1/3,0.1)\approx 0.00065$ is just above $\g_{0}$,
so for the area of the corresponding basin of attraction to come
close to the maximum possible value one needs large values of $\Delta$;
on the contrary $\g(1,0.1)\approx 0.00146$ is not too close to $\g_{0}$ and hence
the area of the corresponding basin of attraction comes closer to the
maximum possible value for smaller $\Delta\approx 10$.

\begin{figure}[ht] 
\centering 
\ifthenelse{\boolean{pdflat}}
{\subfigure {\includegraphics*[width=2.8in]{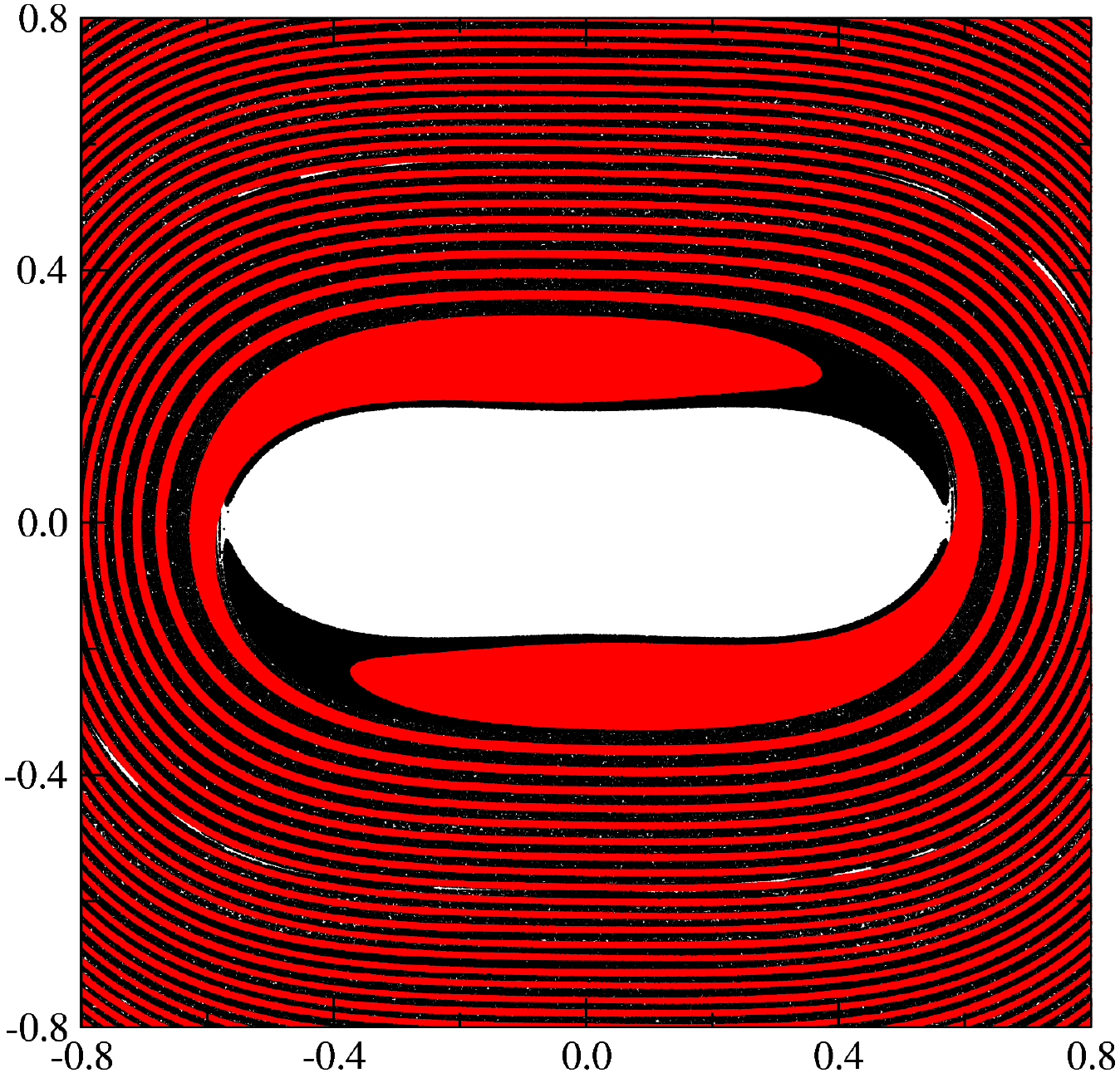} }
\hskip2truecm
\subfigure {\includegraphics*[width=2.8in]{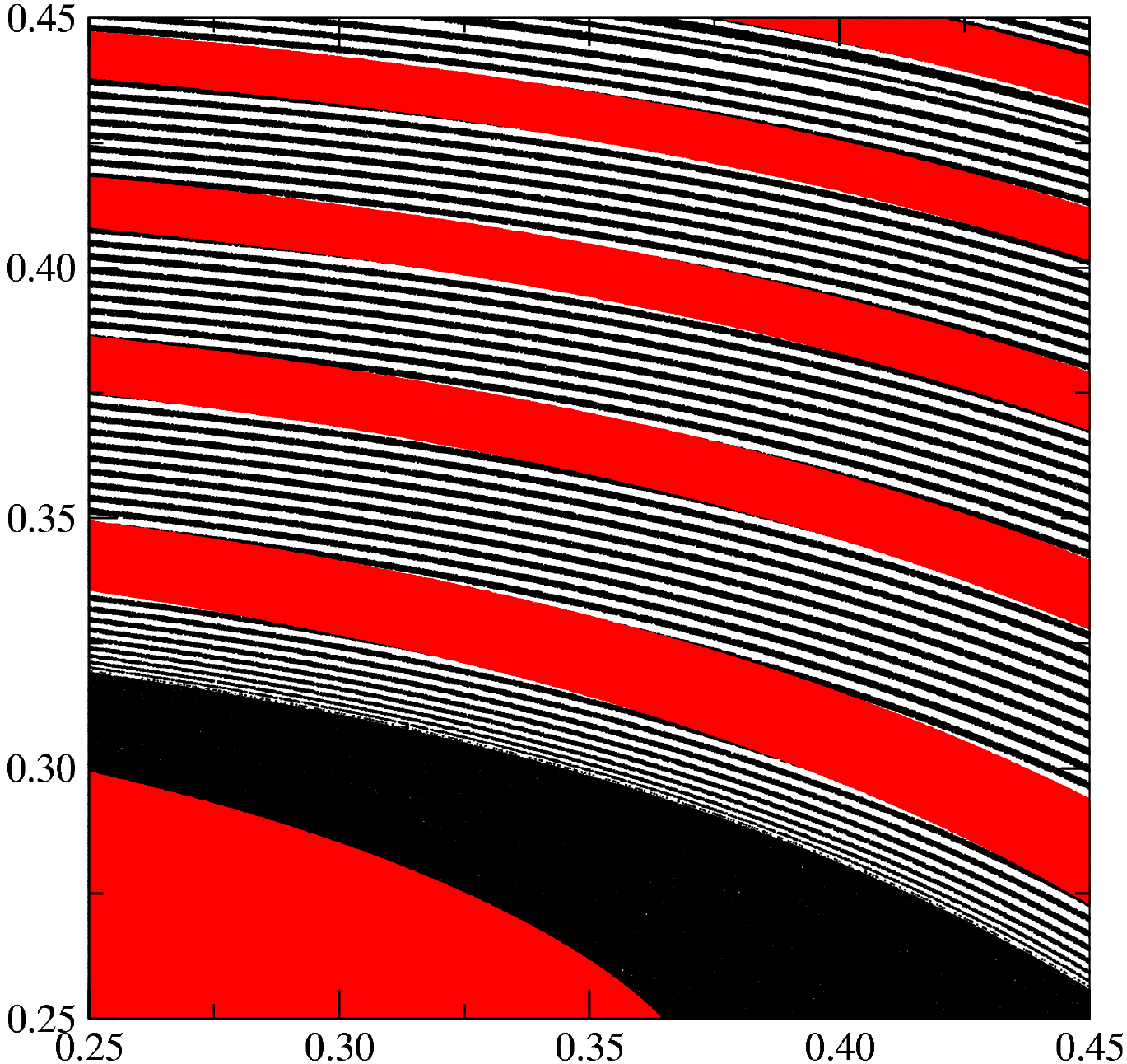} }}
{}  
\caption{Basins of attraction for the 1:2 resonance for constant $\gamma$
(red) and time-varying $\gamma$ (black plus most of the red region). Initial conditions in the white
region either go to the origin or to the 1:4 resonance. The right-hand figure 
shows a portion of the left-hand figure, magnified.}
\label{fig:3.4} 
\end{figure}

We finish this section with a pair of figures showing the difference
between the basins of attraction of the 1:2 resonance for $\gamma =0.005$,
in the cases of constant and time-varying $\g$. According to Table \ref{tab:3.2}, 
changing $\Delta$ from zero --- i.e. constant $\g$ --- to 40 increases the area by
about 8\%, and Figure \ref{fig:3.4} shows how this extra area is distributed
(most of the points of the basin of attraction of the resonance for
constant $\g$ still belong to the basin of attraction for varying $\g$).

\zerarcounters 
\section{The spin-orbit model}
\label{sec:4} 

The spin-orbit model  describes the motion of an asymmetric
ellipsoidal celestial body (satellite)
which  moves in a Keplerian elliptic orbit around a central body
(primary) and rotates around an axis orthogonal to the orbit plane \cite{GP,MD}.
If $\theta$ denotes the angle between the longest axis of the satellite
and the perihelion line,
in the presence of tidal friction the model is described by the equation
\begin{equation} \label{eq:4.1}
\ddot\theta + \e \, G(\theta,t) + \g \, ( \dot \theta -1 ) = 0 ,
\end{equation}
where  $\e,\g>0$ and $\theta\in\TTT=\RRR/2\pi\ZZZ$,
so that the phase space is $\TTT\times\RRR$ (note that \eqref{eq:4.1}
is of the form \eqref{eq:1.2}, with $x=\theta-t$). Here $\e$ is a small parameter,
related to the asymmetry of the equatorial moments of inertia of the satellite,
and $G(\theta,t)=\partial_{\theta} g(\theta,t)$, where
\begin{eqnarray}
g(\theta,t)
& \!\! = \!\! &
\Big( \frac{1}{4}e - \frac{1}{32} e^{3} + \frac{5}{768}e^{5} \Big) \cos (2\theta-t) +
\Big( \frac{1}{2} - \frac{5}{4} e^{2} + \frac{13}{32} e^{4} \Big) \cos (2\theta-2t)
\nonumber \\
& \!\! + \!\! &
\Big( - \frac{7}{4} e + \frac{123}{32}e^{3} - \frac{489}{256} e^{5}\Big) \cos (2\theta-3t) +
\Big( \frac{17}{4}e^{2} - \frac{115}{12} e^{4} \Big) \cos (2\theta-4t) 
\nonumber \\
& \!\! + \!\! &
\Big( - \frac{845}{96} e^{3} + \frac{32525}{1536} e^{5} \Big) \cos (2\theta-5t)  +
\Big( \frac{533}{32} e^{4} \Big) \cos (2\theta-6t)
\label{eq:4.2} \\
& \!\! + \!\! &
\Big( - \frac{228347}{7680} e^{5} \Big) \cos (2\theta-7t) +
\Big( - \frac{1}{96}e^{3} - \frac{11}{1536} e^{5} \Big) \cos (2\theta + t)
\nonumber \\
& \!\! + \!\! &
\Big( \frac{1}{48}e^{4} \Big) \cos (2\theta + 2t) +
\Big( - \frac{81}{2560}e^{5} \Big) \cos (2\theta + 3t) ,
\nonumber
\end{eqnarray}
with $e$ being the eccentricity of the orbit; terms of order $O(e^6)$ have been neglected.
In the celestial mechanics cases  the model \eqref{eq:4.1} may appear oversimplified and more
realistic pictures could be devised \cite{LN,CL1,CL2}. Nevertheless, because
of its simplicity, it is suitable also for analytical investigations (as opposed to 
just numerical ones) on the relevance of friction in the early stages of evolution of 
celestial bodies and for the selection of structurally stable periodic motions.

Values of $e$, $\e$ and $\g$ for some primary-satellite systems of the solar system are
given in Table \ref{tab:4.1}.  The values of $e$ can be found in the literature \cite{MD},
while the derivation of $\e$ and $\g$ is discussed in Appendix \ref{app:b}; see also below.
All satellites of the solar system are trapped in the 1:1 resonance
(rotation period equal to the revolution period), with the remarkable
exception of Mercury (which can be considered as a satellite of the Sun),
which turns out to be in a 3:2 resonance.
The spin-orbit model has been used since the seminal paper by Goldreich and Peale
\cite{GP} in an effort to explain the anomalous behaviour of Mercury.
The ultimate reason is speculated to be related to the large value of the eccentricity.
However, even though higher than for the other primary-satellite systems,
the probability of capture into the 3:2 resonance is still found to be rather low
\cite{GP,CC}. In the following part of this section we aim to investigate
what happens if we take into account the fact that friction increased
during the evolution history of the satellites.

\begin{table}[ht]
\centering
\caption{Values of the constants $e$, $\e$ and $\g$ for some
cases of physical interest for the spin-orbit model \eqref{eq:4.1}.}
\begin{center}
\setlength\tabcolsep{5pt}
\vskip-.5truecm
\vrule
\begin{tabular}{lllllllllllllll}
\hline\noalign{\smallskip}\hline
\null & \vrule\vrule\vrule &
Primary & \vrule\vrule\vrule & Satellite & \vrule\vrule\vrule & $e$ &
\vrule & $\e$ & \vrule & $\g$ \\
\hline\noalign{\smallskip}\hline
E-M & \vrule\vrule\vrule & 
Earth & \vrule\vrule\vrule & Moon & \vrule\vrule\vrule & $\st0.0549$
& \vrule & $\st6.75 \times 10^{-7} $ & \vrule & $\st3.75 \times 10^{-8}$ \\
\hline
S-M & \vrule\vrule\vrule &  
Sun & \vrule\vrule\vrule & Mercury & \vrule\vrule\vrule & $\st0.2056$
& \vrule & $\st8.11 \times 10^{-7} $ & \vrule & $\st3.24 \times 10^{-8}$ \\
\hline
J-G & \vrule\vrule\vrule &  
Jupiter & \vrule\vrule\vrule & Ganymede & \vrule\vrule\vrule & $\st0.0013$ 
& \vrule & $\st4.30\times 10^{-4} $ & \vrule & $\st1.91 \times 10^{-5}$ \\
\hline
J-I & \vrule\vrule\vrule &  
Jupiter & \vrule\vrule\vrule & Io & \vrule\vrule\vrule & $\st0.0041$ 
& \vrule & $\st3.85 \times 10^{-3} $ & \vrule & $\st1.71 \times 10^{-4}$ \\
\hline
S-E & \vrule\vrule\vrule &  
Saturn & \vrule\vrule\vrule & Enceladus & \vrule\vrule\vrule & $\st0.0047$ 
& \vrule & $\st1.41 \times 10^{-2} $ & \vrule & $\st6.26 \times 10^{-4}$ \\
\hline
S-D & \vrule\vrule\vrule &  
Saturn & \vrule\vrule\vrule & Dione & \vrule\vrule\vrule & $\st0.0022$ 
& \vrule & $\st3.85 \times 10^{-3} $ & \vrule & $\st1.71 \times 10^{-4}$ \\
\hline
\end{tabular}
\hspace{-0.1cm}\vrule
\vskip-.5truecm
\end{center}
\label{tab:4.1}
\end{table}

Again one can compute the threshold values of the primary resonances,
by writing $\g(\om)=C_{0}(\om)\e$, up to higher order corrections. If one writes \eqref{eq:4.2} as
\begin{equation} \label{eq:4.3}
g(\theta,t) = \sum_{k\in\ZZZ} a_{k} \cos(2\theta-kt) ,
\end{equation}
one finds  (see Appendix \ref{app:d})
\begin{equation} \label{eq:4.4}
C_{0}(p/q) = \frac{2 q |a_{2p/q}|}{|p-q|} , \qquad \frac{p}{q} \in
\left\{ -1, \pm\frac{1}{2}  , \pm \frac{3}{2}, 2,
\frac{5}{2}, 4 , \frac{7}{2} \right\} ,
\end{equation}
while $C_{0}(1)=\infty$ (that is no threshold value exists for the 1:1 resonance) and $C_{0}(\om)=0$
for any other $\om$; other resonances may appear only at higher order in $\e$.
This leads to the values listed in Table \ref{tab:4.2},
for the primary resonances of the systems considered in Table \ref{tab:4.1}.
Note that for $\g$ large enough, all attractors disappear except
the 1:1 resonance, which becomes a global attractor.

\begin{table}[ht]
\centering
\caption{Values of the constants $C_{0}(p/q)$ for some
primary resonances of the the spin-orbit model \eqref{eq:4.1};
the threshold values are of the form $\g(\om,\e)=C_{0}(\om)\e$.
Only positive $\om$ have been explicitly considered.} 
\begin{center}
\setlength\tabcolsep{5pt}
\vskip-.5truecm
\vrule
\begin{tabular}{llllllllllllllllll}
\hline\noalign{\smallskip}\hline
\null & \vrule\vrule\vrule & $\om$ & \vrule\vrule\vrule &
1/2 & \vrule & 3/2 & \vrule & 2 & \vrule & 5/2 & \vrule & 3 & \vrule & 7/2 \\        
\hline\noalign{\smallskip}\hline
E-M & \vrule\vrule\vrule & $\st C_{0}(\om)$ & \vrule\vrule\vrule
& $\st5.488 \times 10^{-2}$ & \vrule & $\st3.818 \times 10^{-1}$ & \vrule 
& $\st2.545 \times 10^{-2}$ & \vrule & $\st1.928 \times 10^{-3}$ & \vrule
& $\st1.513 \times 10^{-4} $ & \vrule & $\st1.186 \times 10^{-5} $ \\
\hline
S-M & \vrule\vrule\vrule &  $\st C_{0}(\om)$ & \vrule\vrule\vrule
& $\st2.045 \times 10^{-1}$ & \vrule & $\st1.308$ & \vrule
& $\st3.251 \times 10^{-1}$ & \vrule & $\st9.163 \times 10^{-2}$ & \vrule
& $\st2.976 \times 10^{-2}$ & \vrule & $\st8.739 \times 10^{-3}$ \\
\hline
J-G & \vrule\vrule\vrule &  $\st C_{0}(\om)$ & \vrule\vrule\vrule
& $\st1.300 \times 10^{-3}$ & \vrule &  $\st9.100 \times 10^{-3}$ & \vrule
& $\st1.436 \times 10^{-5}$ & \vrule & $\st2.578 \times 10^{-8}$ & \vrule 
& $\st4.757 \times 10^{-11}$ & \vrule & $\st8.832 \times 10^{-14}$ \\
\hline
J-I & \vrule\vrule\vrule &  $\st C_{0}(\om)$ & \vrule\vrule\vrule
& $\st4.100 \times 10^{-3}$ & \vrule &  $\st2.870 \times 10^{-2}$ & \vrule 
& $\st1.429 \times 10^{-4}$ & \vrule & $\st8.088 \times 10^{-7}$ & \vrule 
& $\st4.707 \times 10^{-9}$ & \vrule & $\st2.756 \times 10^{-11}$ \\
\hline
S-E & \vrule\vrule\vrule &  $\st C_{0}(\om)$ & \vrule\vrule\vrule
& $\st4.700 \times 10^{-3}$ & \vrule &  $\st3.290 \times 10^{-2}$ &\vrule 
& $\st1.878 \times 10^{-4}$ & \vrule & $\st1.218 \times 10^{-6}$ & \vrule 
& $\st8.128 \times 10^{-9}$ & \vrule & $\st5.455 \times 10^{-11}$ \\
\hline
S-D & \vrule\vrule\vrule &  $\st C_{0}(\om)$ & \vrule\vrule\vrule
& $\st2.200 \times 10^{-3}$ & \vrule & $\st1.540 \times 10^{-2}$ & \vrule
& $\st4.114 \times 10^{-5}$ & \vrule & $\st1.259\times 10^{-7}$ & \vrule 
& $\st3.902 \times 10^{-10}$ & \vrule & $\st1.226 \times 10^{-12} $ \\
\hline
\end{tabular}
\hspace{-0.1cm}\vrule
\vskip-.5truecm
\end{center}
\label{tab:4.2}
\end{table}

It seems reasonable on physical grounds (see Appendix \ref{app:e})
to assume that friction was increasing in the past up to the present-day value
$\g_{0}=\g$, with $\g$ as in Table \ref{tab:4.1}.
Application to the spin-orbit model for S-M would require numerics with
very small values of $\e$ and $\g_{0}$: the discussion in Appendix \ref{app:b}
provides the values in Table \ref{tab:4.1}, so that $\g_{0} \sim 0.05\,\e$.
However, the value usually taken for $\e$ in the literature is $\e \sim 10^{-4}$
(see Appendix \ref{app:b}). In both cases, $\g_{0}$ is far below the
threshold value $\g(\om,\e)$, especially
for the most interesting resonance $\om=3/2$ \cite{GP,MD},
as $\g(3/2,\e) \sim 1.3 \, \e$ (see Table \ref{tab:4.2});
we refer to Section \ref{sec:6} for further comments.

As already noted in \cite{CC}, the small value of $\g_{0}$ represents a serious difficulty
from a numerical point of view, because it requires very long integration
(for a very large number of initial conditions --- see comments in Appendix \ref{app:f}).
Nevertheless, the discussion in Section \ref{sec:3} about the driven cubic oscillator allows us 
to draw the following conclusions about the spin-orbit model.

The results in Table \ref{tab:2.3} suggest that the relative areas $A(\om,\g,\e)$ of the
basins of attraction are almost constant for values of $\g$ much smaller
than the threshold values; more precisely they assume values close to
$A_{\rm max}(\om,\e)$. In the case of increasing friction, if the final value $\g_{0}$
is much smaller than the threshold value, the basin turns out to
have more or less the same size close to $A_{\rm max}(\om,\e)$ as it would have if $\g$ were set
equal to $\g_{0}$ since the beginning. The same scenario is expected
in the case of the spin-orbit problem. In particular a crucial aspect is understanding
when the final value can be considered `much smaller' than the
threshold value or comparable to it. Again the analysis in Section \ref{sec:3} is useful:
pragmatically, we shall define $\g$ much smaller than $\g(\om,\e)$
when $A(\om,\g,\e)$ is close to the maximum possible value $A_{\rm max}(\om,\e)$. 
Therefore it becomes fundamental to check whether, for the current values of $\g$ and $\e$
as given in the literature, $A(\om,\g,\e)$ is either close to $A_{\rm max}(\om,\e)$ or much smaller.

\begin{enumerate}

\item
If we assumed $A(\om,\g,\e)$ to be close to $A_{\rm max}(\om,\e)$
the time-dependence of friction would not change the general picture as observed today.
From this point of view, our results would be a bit disappointing: indeed
the relative area of the basin of attraction of the 3:2 resonance,
with the values of $\e$ and $\g$ usually taken in the literature, is found to be
rather small for S-M \cite{CC} and including the time-dependence of friction
in the analysis would not give larger estimates.
On the other hand, also in this second case, our analysis 
would provide some more information: it would yield
that the results available in the literature \cite{CL1,CL2,CC}
would remain correct even if time-dependent friction were included. 
In particular, to explain why Mercury has been captured into the
3:2 resonance, other mechanisms should be be invoked,
such as the chaotic evolution of its orbit \cite{CL1}.
Of course the values of $\e$ and $\g$ used in the literature
are only speculative: again our analysis suggests that
the results would not change in a sensible way even by taking
different values for one or both parameters --- see also comments in Section \ref{sec:6}.

\item
On the contrary if $A(\om,\g,\e)$ were much smaller than $A_{\rm max}(\om,\e)$,
taking into account the time-dependence of
friction would imply a larger basin of attraction with respect to the case of constant friction. 
In this case the the exact values of the parameters $\e$ and $\g$
would play a fundamental role --- again see also Section \ref{sec:6}.

\end{enumerate}

\zerarcounters 
\section{Remarks and comments}
\label{sec:5} 

\subsection{Different values of the perturbation parameter}
In Sections \ref{sec:2} and \ref{sec:3} we have fixed $\e=0.1$. However,
by taking different values of $\e$, the phenomenology does not change.
For instance, for $\e=0.5$ and $\e=0.01$ we have the relative areas
listed in Tables \ref{tab:5.1} and \ref{tab:5.2} (for $\e=0.5$ only a
few attractors are taken into account in Table \ref{tab:5.1}).

\begin{table}[ht]
\centering
\caption{Relative areas of the parts of the basins of attraction
contained inside the square $\calQ$ for $\e=0.5$.
($\om=0$ denotes the origin). $1\,000\,000$ random initial conditions have been taken in $\calQ$.}
\begin{center}
\setlength\tabcolsep{5pt}
\vskip-.5truecm
\vrule
\begin{tabular}{lllllllllllllll}
\hline\noalign{\smallskip}\hline
$\om$ & \vrule\vrule\vrule & $0$ & \vrule & $1/2$ & \vrule & $1/4$ &
\vrule & $1$a & \vrule & $1$b & \vrule & $1/6$ \\
\hline\noalign{\smallskip}\hline
$\g = 0.1000$ & \vrule\vrule\vrule & \!\!\!100.0 & \vrule & 00.0 & \vrule & 00.0
& \vrule & 00.0 & \vrule & 00.0 & \vrule & 00.0 \\
\hline
$\g = 0.0750$ & \vrule\vrule\vrule & 70.5 & \vrule & 29.5 & \vrule & 00.0 
& \vrule & 00.0  & \vrule & 00.0 & \vrule & 00.0 \\
\hline
$\g = 0.0500$ & \vrule\vrule\vrule & 49.8 & \vrule & 50.2 & \vrule & 00.0 
& \vrule & 00.0 & \vrule & 00.0 & \vrule & 00.0 \\
\hline
$\g = 0.0250$ & \vrule\vrule\vrule & 32.0 & \vrule & 56.0 & \vrule & 05.8
& \vrule & 03.1 & \vrule & 03.1 & \vrule & 00.0 \\
\hline
$\g = 0.0050$ & \vrule\vrule\vrule & 10.4 & \vrule & 48.8 & \vrule & 07.7
& \vrule & 10.5 & \vrule & 10.5 & \vrule & 00.0 \\
\hline
$\g = 0.0025$ & \vrule\vrule\vrule & 08.1 & \vrule & 36.3 & \vrule & 06.8
& \vrule & 11.5 & \vrule & 11.5 & \vrule & 00.8 \\
\hline
$\g = 0.0010$ & \vrule\vrule\vrule & 06.9 & \vrule & 37.5 & \vrule & 04.4
& \vrule & 11.4 & \vrule & 11.4 & \vrule & 01.6 \\
\hline
\end{tabular}
\hspace{-0.1cm}\vrule
\vskip-.5truecm
\end{center}
\label{tab:5.1}
\end{table}

\begin{table}[ht]
\centering
\caption{Relative areas of the parts of the basins of attraction
contained inside the square $\calQ$ for $\e=0.01$.
($\om=0$ denotes the origin). $500\,000$ random initial conditions have been taken in $\calQ$.}
\begin{center}
\setlength\tabcolsep{5pt}
\vskip-.5truecm
\vrule
\begin{tabular}{lllllllllllllll}
\hline\noalign{\smallskip}\hline
$\om$ & \vrule\vrule\vrule & $0$ & \vrule & $1/2$ & \vrule & $1/4$ &
\vrule & $1$a & \vrule & $1$b & \vrule & $1/6$ \\
\hline\noalign{\smallskip}\hline
$\g = 0.0020$ & \vrule\vrule\vrule & \!\!\!100.0 & \vrule & 00.0 & \vrule & 00.0
& \vrule & 00.0 & \vrule & 00.0 & \vrule & 00.0 \\
\hline
$\g = 0.0015$ & \vrule\vrule\vrule & 98.1 & \vrule & 01.9 & \vrule & 00.0 
& \vrule & 00.0 & \vrule & 00.0 & \vrule & 00.0 \\
\hline
$\g = 0.0010$ & \vrule\vrule\vrule & 94.0 & \vrule & 06.0 & \vrule & 00.0 
& \vrule & 00.0 & \vrule & 00.0 & \vrule & 00.0 \\
\hline
$\g = 0.0005$ & \vrule\vrule\vrule & 88.3 & \vrule & 10.6 & \vrule & 01.1 
& \vrule & 0.00 & \vrule & 0.00 & \vrule & 00.0 \\
\hline
$\g = 0.0001$ & \vrule\vrule\vrule & 78.2 & \vrule & 15.3 & \vrule & 06.5
& \vrule & 00.0 & \vrule & 00.0 & \vrule & 00.0 \\
\hline
\end{tabular}
\hspace{-0.1cm}\vrule
\vskip-.5truecm
\end{center}
\label{tab:5.2}
\end{table}

The general scenario is the same as in the case $\e=0.1$,
with obvious quantitative differences due to the fact that
for smaller values of $\e$ (say $\e=0.01$) only primary resonances are relevant 
unless $\g$ is very small, while for larger $\e$
(say $\e=0.5$) more and more resonances appear by taking smaller values of $\g$
(because powers of $\e$ are not much smaller than $\e$ itself).
Indeed, if $\e=0.5$, even for $\g=0.005$ we detect numerically more than 50 attractors and the
classification of resonance $\om$ according to the the value of $n(\om)$
becomes meaningless. For instance, for $\e=0.5$, one has  $A(3/4,0.0025,0.5)\approx 11.4\%$
and $A(3/4,0.001,0.5) \approx 6.8\%$, to be compared with the values
for $\om=1/4$ in Table \ref{tab:5.1}.
Moreover for large values of $\e$, say for $\e=0.5$,
the bending of the curves $A(\g,\om,\e)$ in Figure \ref{fig:2.2} is more pronounced
and the monotonic decrease observed for $\e=0.1$ when $\g$ tends to $0$
seems to be violated (compare the values $A(1/2,\g,0.5)$ for $\g=0.0025$ 
and $\g=0.001$ in Table \ref{tab:5.1}); see also the comments in Section \ref{sec:6}.

\subsection{Different functions $\boldsymbol{\g}\boldsymbol{(}\boldsymbol{t}\boldsymbol{)}$}

\begin{figure}[ht] 
\centering 
\ifthenelse{\boolean{pdflat}}
{\subfigure{\includegraphics*[width=2.8in]{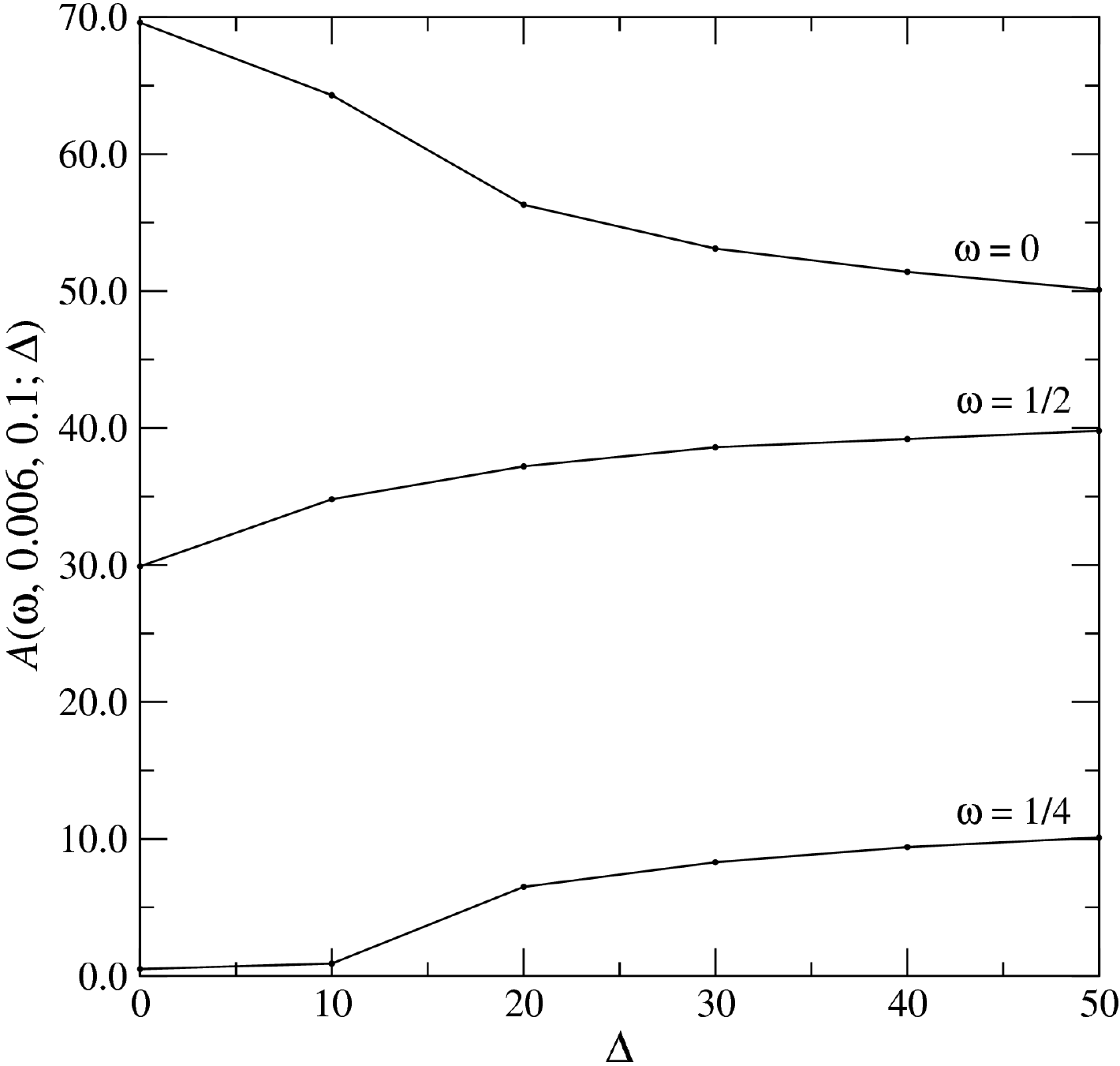}
\hskip2.truecm
\subfigure{\includegraphics*[width=2.8in]{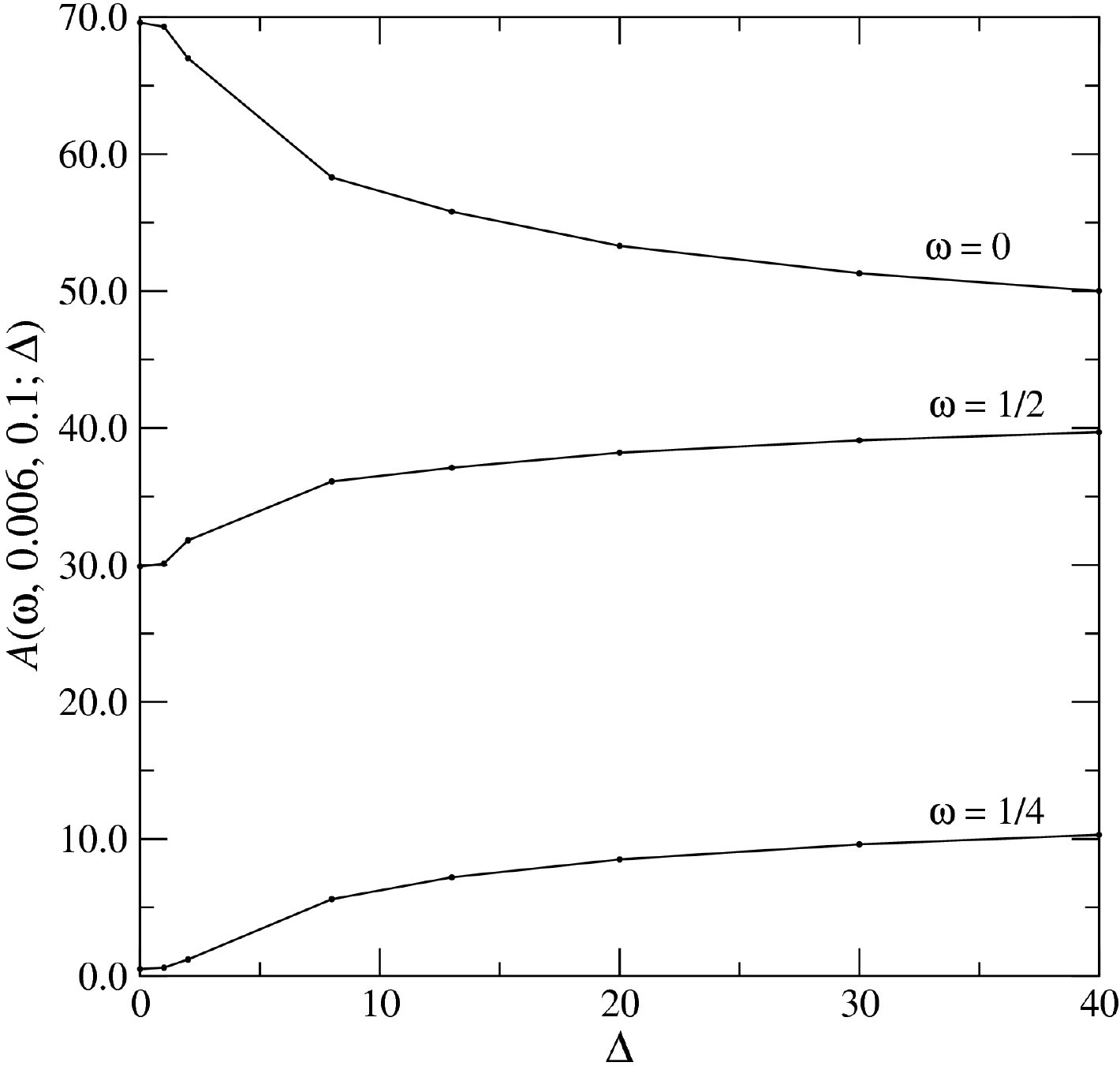} }}}
{\subfigure{\includegraphics*[width=2.8in]{fig51.eps}
\hskip2.truecm
\subfigure{\includegraphics*[width=2.8in]{fig52.eps} }}}
\caption{Relative areas of the basins of attraction versus $\Delta$
for $\g(t)$ given by (\ref{eq:3.2}), with $\g_{0}=0.006$ (left-hand figure)
and for $\g(t)$ given by (\ref{eq:5.1}), with $\g_{0}=0.006$ (right-hand figure).}
\label{fig:5.1} 
\end{figure}

As stated at the beginning of Section \ref{sec:3}, the exact form of
the function $\g(t)$ should not be relevant. As a check we studied (\ref{eq:3.1})
with $\g(t)$ given by both (\ref{eq:3.2}) and
\begin{equation}
\gamma(t)  = \gamma_{0} \big( 1 - {\rm e}^{-t/T_{0}} \big) ,
\label{eq:5.1}
\end{equation}
where $\gamma_{0}$ and $T_{0}$ are positive constants,
by setting $T_{0}=\Delta/\g_{0}$ and changing $\Delta$ for fixed values of $\g_{0}$.
The results show that the same behaviour is obtained in both cases.
For instance, for $\g_{0}=0.006$, one has the results in Tables \ref{tab:5.3} and \ref{tab:5.4};
see also Figure \ref{fig:5.1}.

\begin{table}[ht]
\centering
\caption{Relative areas $A(\om,0.006,0.1;\Delta)$ of the parts of the basins of attraction
contained inside the square $\calQ$
for $\e=0.1$ and $\g(t)$ given by (\ref{eq:3.2}) with $\g_{0}=0.006$
and $T_{0}=\Delta/\g_{0}$, for various values of $\Delta$ and $\om=0,1/2,1/4$.
($\om=0$ denotes the origin). $500\,000$ random initial conditions have been taken in $\calQ$.}
\begin{center}
\setlength\tabcolsep{5pt}
\vskip-.5truecm
\vrule
\begin{tabular}{lllllllllllllllllllllllll}
\hline\noalign{\smallskip}\hline
$\Delta$ & \vrule\vrule\vrule & $0$ & \vrule & $10$ & \vrule 
& $20$ & \vrule & $30$ & \vrule & $40$ & \vrule & $50$ \\ 
\hline\noalign{\smallskip}\hline
   $\om = 0$ & \vrule\vrule\vrule & $69.6$ & \vrule & $64.3$ & \vrule 
& $56.3$ & \vrule & $53.1$ & \vrule & $51.4$ & \vrule & $50.1$ \\ 
\hline
$\om = 1/2$ & \vrule\vrule\vrule & $29.9$ & \vrule & $34.8$ & \vrule
& $37.2$ & \vrule & $38.6$ & \vrule & $39.2$ & \vrule & $39.8$ \\ 
\hline
$\om = 1/4$ & \vrule\vrule\vrule & $00.5$ & \vrule & $00.9$ & \vrule 
& $06.5$ & \vrule & $08.3$ & \vrule & $09.4$ & \vrule & $10.1$ \\ 
\hline
\end{tabular}
\hspace{-0.1cm}\vrule
\vskip-.5truecm
\end{center}
\label{tab:5.3}
\end{table}

\begin{table}[ht]
\centering
\caption{Relative areas $A(\om,0.006,0.1;\Delta)$ of the parts of the basins of attraction
contained inside the square $\calQ$
for $\e=0.1$ and $\g(t)$ given by (\ref{eq:5.1}) with $\g_{0}=0.006$
and $T_{0}=\Delta/\g_{0}$, for various values of $\Delta$ and $\om=0,1/2,1/4$
($\om=0$ denotes the origin). $500\,000$ random initial conditions have been taken in $\calQ$.}
\begin{center}
\setlength\tabcolsep{5pt}
\vskip-.5truecm
\vrule
\begin{tabular}{lllllllllllllllllllllll}
\hline\noalign{\smallskip}\hline
$\Delta$ & \vrule\vrule\vrule & $0$ & \vrule & $1$ & \vrule & $2$ & \vrule & $8$
& \vrule & $13$ & \vrule & $20$ & \vrule & $30$ & \vrule & $40$ \\
\hline\noalign{\smallskip}\hline
   $\om = 0$ & \vrule\vrule\vrule & $69.6$ & \vrule & $69.3$ & \vrule & $67.0$ & \vrule & $58.3$ 
& \vrule & $55.8$ & \vrule & $53.3$ & \vrule & $51.3$ & \vrule & $50.0$ \\ 
\hline
$\om = 1/2$ & \vrule\vrule\vrule & $29.9$ & \vrule & $30.1$ & \vrule & $31.8$ & \vrule & $36.1$ 
& \vrule & $37.1$ & \vrule & $38.2$ & \vrule & $39.1$ & \vrule & $39.7$ \\ 
\hline
$\om = 1/4$ & \vrule\vrule\vrule & $00.5$ & \vrule & $00.6$ & \vrule & $01.2$ & \vrule & $05.6$ 
& \vrule & $07.2$ & \vrule & $08.5$ & \vrule & $09.6$ & \vrule & $10.3$ \\ 
\hline
\end{tabular}
\hspace{-0.1cm}\vrule
\vskip-.5truecm
\end{center}
\label{tab:5.4}
\end{table}

In particular, in both cases the relative areas $A(\om,\g_{0},0.1;\Delta)$ start at the values
corresponding to constant dissipation $\g=\g_{0}$ for $\Delta=0$ and
then either decrease (for the origin) or increase (for the periodic orbits),
apparently to some asymptotic value close to $A_{\rm max}(\om,0.1)$.
For instance the relative areas corresponding to (\ref{eq:3.2})
with $\Delta=30$ are very close to those corresponding to (\ref{eq:5.1})
with $\Delta=20$ (compare Table \ref{tab:5.3} with Table \ref{tab:5.4}).
The asymptotic values seem to be the same in both cases.

\zerarcounters 
\section{Conclusions and open problems}
\label{sec:6} 

In this paper we have studied how the slow growth of friction may
affect the asymptotic behaviour of dissipative dynamical systems.
We have focused on a simple paradigmatic model, the periodically driven
cubic oscillator, particularly suited for numerical investigations.
Nevertheless we think that the results hold in the more general
setting considered in Section \ref{sec:1}.
The main result, discussed in Section \ref{sec:3}, can be summarised as follows:
on the one hand it is the final value of the damping coefficient that
determines which attractors are present, but on the other hand the sizes
of the corresponding basins of attraction strongly depend
on the full evolution of the damping coefficient itself, in particular on its growth rate.
Let $\g(p/q,\e)$ and $\g_{0}$ denote the threshold value of the $p:q$ resonance 
and the final value of the damping coefficient, respectively.
If $\g_{0}>\g(p/q,\e)$ the attractor disappears. Otherwise the following possibilities arise: 
if $\g_{0} \ll \g(p/q,\e)$ the area of the corresponding basin of attraction is
more or less the same as in the case with constant damping coefficient equal to $\g_{0}$,
whereas it is larger if $\g_{0} \lesssim \g(p/q,\e)$.
In the latter case, the closer $\g_{0}$ is to  $\g(p/q,\e)$, the slower the growth rate
of $\g(t)$ required for the maximum possible area of the basin of attraction to be attained.
Moreover when $\g_{0} \ll \g(p/q,\e)$ the area can be even a little
smaller than what it would be if the damping were fixed at $\g_{0}$
since $t = 0$, because other attractors may have acquired a larger
basin of attraction while the damping coefficient increased.
It is not possible for the relative area to be larger than the value $A_{\rm max}(p/q,\e)$,
which therefore represents an upper bound.

Finally, let us mention a few open problems which would deserve
further investigation, also in relation with the spin-orbit problem.

\begin{enumerate}

\item As far as a dynamical system can be considered as a perturbation of an integrable one,
all attractors seem to be either equilibrium points or periodic orbits:
at least, this is what emerges from numerical simulations. It would be
interesting to have a proof of this behaviour, even for some simple
model such as (\ref{eq:2.1}), in particular of the fact
that neither strange attractors nor periodic solutions other than
subharmonic ones appear.

\item The analysis presented in Section \ref{sec:3} covers small values of $\g$,
but still not as small as desirable. It would be worthwhile to study
the behaviour of the curves $A(\om,\g,\e)$ in the limit of even smaller $\g$,
for instance by implementing some numerical integrator which allows
us to decrease the running time of the programs without losing accuracy in the results.

\item In studying the spin-orbit model in Section \ref{sec:4}
if one really wanted to consider the past history of the system,
then not only $\g$ but also $\e$ should be taken to depend on time:
this would introduce further difficulties and a more detailed
model for the evolution of the satellite would be needed
(see also comments at the end of Appendix \ref{app:b}).
Of course, this would be by no means an easy task.
The use of a model for the time-dependence of friction already
raises several problems, as we highlight in Appendix \ref{app:e}.

\item We have seen that, for the spin-orbit model, in the case of constant friction
the exact value of the parameters $\e$ and $\g$ is fundamental.
For instance if $\g=\g_{1}$ such that $A(\om,\g_{1},\e) \approx A_{\rm max}(\om,\e)$,
then we have a basin of attraction much larger than for $\g=\g_{2}$, where
$\g_{2} \gg \g_{1}$ is close to the threshold value $\g(\om,\e)$.
On the contrary, our analysis shows that, when one takes into account
that friction slowly increased during the solidification process of the satellite,
then for both values $\g_{1}$ and $\g_{2}$ we expect more or less
the same relative area close to $A_{\rm max}(\om,\e)$.
This is useful information because it shows that
exact values of the parameters $\e$ and $\g$ are not 
essential in the case of time-dependent friction, as long
as $\g$ is not much smaller than the threshold value
(we mean `much smaller' in the sense of Section \ref{sec:4}):
of course the values of $\g$ and $\e$ are essential
if $\g$ is much smaller than $\g(\om,\e)$ --- see next item.

\item As remarked in Section \ref{sec:4}, to study the probability
of capture of Mercury into the 3:2 resonance, it becomes
crucial to determine the value $A_{\rm max}(3/2,\e)$ and 
the current values of $\g$ and $\e$ to see if
$\g$ turns out to be much smaller than $\g(3/2,\e)$,
that is if one has $A(3/2,\g) \ll A_{\rm max}(3/2,\e)$.
If this were the case, by using time-dependent friction, the relative area
of the basin of attraction of the 3:2 resonance would be much larger than
what was found for constant damping coefficient  $\g$
(estimated around $13\%$ in \cite{CC}).
By assuming the values of $\g$ given in the literature, for
which $\g\sim 10^{-8}$ and $\e \sim 10^{-4}$, and using
that $\g(3/2,\e) \sim 1.3 \, \e$, the relation $A(3/2,\g) \ll A_{\rm max}(3/2,\e)$
could be verified only by assuming for $A(3/2,\g,\e)$ a very slow variation in $\g$
for $\e \sim 10^{-4}$. This does not seem impossible:
already for the driven cubic oscillator the profiles of the curves $A(\om,\g,\e)$
seem to have a much smoother variation for smaller values of $\e$
(compare Table \ref{tab:2.3} for $\e=0.1$ with Table \ref{tab:5.1}
for $\e=0.5$), so it could happen then by decreasing $\e$ further
the curves $A(\om,\g,\e)$ could be nearly flat on much longer
intervals: in other words, by fixing $\e\sim10^{-4}$ and decreasing $\g$,
the curve $A(\om,\g,\e)$ could have not yet reached its maximum value
$A_{\max}(\om,\e)$ at $\g\sim10^{-8}$.
As a consequence, the exact values of the parameters $\e$ and $\g$
and the profile of the curve $A(3/2,\e,\g)$ could be fundamental.
In any case, it would be very interesting to study numerically the
spin-orbit model for very small values of $\e$ and $\g$,
in order to obtain the profiles of the curves $A(3/2,\g,\e)$ versus $\g$.

\end{enumerate}

\noindent \textit{Acknowledgments.} We are grateful to Giovanni Gallavotti
for very useful discussions and critical remarks, especially
on the spin-orbit model and the formation and evolution of celestial bodies.

\appendix

\zerarcounters
\section{Some analytical results on system (\ref{eq:2.1})}
\label{app:a}

\subsection{Global attraction to the origin for large $\boldsymbol{\g}$}

For $\e$ small enough introduce the positive function $F(t)$ such that
$F^{2}(t)=1+\e f(t)$ and rescale time through the
Liouville transformation \cite{MW,BDGG}
\begin{equation}
\tau := \int_{0}^{t} {\rm d}s \, F(s) .
\label{eq:a.1}
\end{equation}
Then we can rewrite (\ref{eq:2.1}) as
\begin{equation}
\begin{cases}
x' = y , & \\
y' = \displaystyle{ - x^{3} - \frac{y}{F(t)} 
\Big( \g + \frac{F'(t)}{F(t)} \Big) } , & \end{cases}
\label{eq:a.2}
\end{equation}
where the prime denotes derivative with respect to $\tau$. Define
\begin{equation}
I(x,\dot x, t) := \frac{y^{2}}{2} + \frac{x^{4}}{4} ,
\label{eq:a.3}
\end{equation}
which is an invariant for (\ref{eq:2.1}) with $\e=\g=0$.
More generally one has
\begin{equation}
I' = - \frac{1}{F(t)} \left( \g + \frac{\e \dot f(t)}{2(1+\e f(t))} \right) y^{2} ,
\nonumber
\end{equation}
so that, if
\begin{equation}
\g > - \min \frac{\e\dot f(t)}{2(1+\e f(t))} ,
\label{eq:a.4}
\end{equation}
by using Barbashin-Krasovsky-La Salle's theorem \cite{K}, we find that the origin is
an asymptotically stable equilibrium point and every initial datum is attracted to it as $t\to\infty$.

\subsection{Rate of convergence to attractors} 
The periodic orbits for the system (\ref{eq:2.1}) appear in pairs of stable and unstable orbits:
this is a consequence of Poincar\'e-Birkhoff's theorem \cite{Birkhoff}. 
Let us consider the primary resonances, so that we can set $\g=\e \, C$,
with $C$ a constant independent of $\e$.

We rewrite the equations of motion in action-angle variables $(I,\f)$ \cite{BBDGG}
\begin{equation}
\begin{cases}
\dot \f = (3I)^{1/3} + \e  (3I)^{1/3} f(t) \, \cn^{4} \f - \e \, C \cn\f \, \sn \f \, \dn \f , & \\
\dot I = \e  (3I)^{4/3} f(t) \, \cn^{3} \f \, \sn \f \, \dn \f - \e \,
C  (3I) \, \sn^{2}\f \, \dn^{2} \f , & 
\end{cases}
\label{eq:a.5}
\end{equation}
where $\cn\f,\sn\f,\dn\f$ are the cosine-amplitude, sine-amplitude, delta-amplitude
functions, respectively, with elliptic modulus $k=1/\sqrt{2}$ \cite{GR}.

Let $K(k)$ be the complete elliptic integral of the first type.
For $\e=0$ the periodic solution to (\ref{eq:2.1}) with frequency $\om=p/q$ is of the form
$x(t) = \al \, \cn(\al(t+t_{0}))$, with $2\pi\al=4\om K(1/2)$ and $t_{0}$ suitably fixed \cite{BBDGG}.
In terms of action-angle variable this gives $I=I_{0}:=\al^3/3$ and $\f=\f_{0}(t):=\al (t+t_{0})$.

Linearisation of (\ref{eq:a.5}) around the periodic solution leads to the system
\begin{equation}
\left( \begin{matrix} \dot {\de\f} \\ \dot {\de I} \end{matrix} \right) =
L(t) \left( \begin{matrix} \de\f \\ \de I \end{matrix} \right) , \qquad
L(t)=L_{0}+ \e L_{1}(t) + O(\e^{2}) , 
\label{eq:a.6} 
\end{equation}
where
\begin{equation}
L_{0}(t) = \left( \begin{matrix} 0 & \al^{-2} \\ 0 & 0 \end{matrix} \right) , \qquad
L_{1}(t) = \left( \begin{matrix} L_{11}(t) & L_{12}(t) \\ L_{21}(t) & L_{22}(t) \end{matrix} \right) ,
\label{eq:a.7}
\end{equation}
with
\begin{eqnarray}
& & L_{11}(t) = f(t) \, \dot a(t) - C \, \al^{-1} \dot b(t) , \qquad
L_{12}(t) = -2\al^{-5} I_{1}(t) + \al ^{-2} f(t) \, a(t) , \nonumber \\
& & L_{21}(t) =  \al^{3} f(t) \, \dot c(t) - C \, \al^{2} \dot d(t) ,\qquad 
L_{22}(t) = 4 \al \,  f(t) \, c(t) - 3 C \, d(t) , 
\label{eq:a.8}
\end{eqnarray}
where we have defined
\begin{eqnarray}
& & a(t) := \cn^{4} \f_{0}(t) ,  \qquad
b(t) := \cn \f_{0}(t) \, \sn \f_{0}(t) \, \dn \f_{0}(t) , \nonumber \\
& & c(t) := \cn^{3} \f_{0}(t) \, \sn \f_{0}(t) \, \dn \f_{0}(t) , \qquad
d(t) := \sn^{2}\f_{0}(t) \, \cn^{2} \f_{0}(t) 
\label{eq:a.9}
\end{eqnarray}
and denoted by $(\f_{1}(t),I_{1}(t))$ the first order of the periodic solution.

Let us denote by $ W(t) = W_{0}(t)+\e W_{1}(t)+O(\e^{2})$ the Wronskian matrix,
that is the matrix whose columns are two independent solutions of the
linearised system (\ref{eq:a.2}), so that $\dot W(t)=L(t)W(t)$, with $W(0)=\uno$. Then one has
\begin{equation}
W_{0}(t) = \exp t L_{0} = \left( \begin{matrix}
1 & \al^{-2} t \\ 0 & 1 \end{matrix} \right) 
\label{eq:a.10}
\end{equation}
while $W_{1}(t)$ is obtained by solving the system $\dot W_{1} = L_{0} W_{1} + L_{1} W_{0}(t) $, i.e.
\begin{equation}
W_{1}(t) = W_{0}(t) \left[
W_{1}(0) + \int_{0}^{t} \der \tau
\left( W_{0}(\tau) \right)^{-1} L_{1}(\tau) \, W_{0}(\tau) \right] ,
\label{eq:a.11}
\end{equation}
where one has to take $W_{1}(0)=0$ in order to have $W(0)=\uno$.

A trivial computation shows that in (\ref{eq:a.10}) one has
\begin{equation}
W_{0}(t) ^{-1} L_{1}(t) W_{0}(t) =
\left( \begin{matrix}
L_{11}(t) - \al^{-2} t \, L_{21}(t) &
L_{12}(t) + \al^{-2} t \big( L_{11}(t) - L_{22}(t) - \al^{-2}  L_{21}(t) \big) \\
L_{21}(t) & L_{22}(t)  + \al^{-2} t \, L_{21}(t) \end{matrix} \right) .
\nonumber
\end{equation}

Let $T=2\pi q$ be the period of the periodic solution.
The Floquet multipliers around the periodic solution
are the eigenvalues of the Wronskian matrix, computed at time $T$.
Denote by $x^{k}(t)$ the $k$th primitive of any function $x(t)$
with $x^{k}(0)=0$ (so that $\dot x^{k}(t) = x^{k-1}(t)$,
with $x^{0}(t)=x(t)$). Then, by using that
\begin{eqnarray}
& & \displaystyle{ \int_{0}^{T} \der t \, x(t) = x^{1}(T) } , \qquad
\displaystyle{ \int_{0}^{T} \der t \, t \, x(t) = Tx^{1}(T) - x^{2}(T) } , \nonumber \\
& & \displaystyle{ \int_{0}^{T} \der t \, t^{2} x(t) = T^{2}x^{1}(T) -
2Tx^{2}(T) + 2 x^{3}(T) } ,
\label{eq:a.12}
\end{eqnarray}
we obtain that
\begin{equation}
W_{1}(T) =
\left( \begin{matrix}
L_{11}^{1} (T) + \al^{-2} L_{21}^{2}(T) & 
L_{12}^{1}(T) + \al^{-2} \big( T L_{11}^{1}(T) - L_{11}^{2} + 
L_{22}^{2} + \al^{-2} (T L_{21}^{2} - 2 L_{21}^{3}(T)) \big) \\
L_{21}^{1}(T) & L_{22}^{1}(T)  + \al^{-2} \, ( T L_{21}^{1}(T) - L_{21}^{2}(T) )
\end{matrix} \right).
\nonumber
\end{equation}
For $\e=0$, the corresponding Floquet multipliers are equal to $1$.
To first order they are the roots $\la_{\pm}$ of the
equation $\la^{2} - 2b_{0} \la + c_{0}=0$, with
\begin{equation}
b_{0} := 1 + \frac{\e}{2} \left( L_{11}^{1}(T) + L_{22}^{1}(T) + \al^{-2} T L_{21}^{1}(T) \right) ,
\quad c_{0} := 1 + \e \left( L_{11}^{1}(T) + L_{22}^{1}(T) \right) ,
\label{eq:a.13}
\end{equation}
so that
\begin{equation}
\la_{\pm} = 1 \pm \sqrt{\e\,\al^{-2} T L_{21}^{1}(T)} +
\frac{\e}{2} \left[ L_{11}^{1}(T) + L_{22}^{1}(T) + \al^{-2} T
L_{21}^{1}(T) \right] + o(\e) .
\label{eq:a.14}
\end{equation}
One has
\begin{equation}
L_{11}^{1}(T)+L_{22}^{1}(T) = \int_{0}^{T} \der t \,
f(t) \big( \dot a(t) + 4 \al \, c(t) \big) - C
\int_{0}^{T} \der t \, \big( \al^{-1} \dot b(t) + 3 d(t) \big)
\label{eq:a.15} 
\end{equation}
One immediately realises that the first integral vanishes and hence
\begin{equation}
L_{11}^{1}(T)+L_{22}^{1}(T) = - 3 C  T \mu , \qquad
\mu := \frac{1}{T} \int_{0}^{T} \der t \, d(t) > 0 ,
\label{eq:a.16}
\end{equation}
while, for the stable periodic solution, $t_{0}$ is such that $L_{21}^{1}(T)<0$.
Therefore the Floquet multipliers are of the form
\begin{equation}
\la_{\pm} = 1 \pm \ii \la_{0} \sqrt{T \e} - \left( \la_{0}^{2} + 3 C \mu \right) T \e + o(\e) ,
\nonumber
\end{equation}
with $\la_{0}>0$. The corresponding Lyapunov exponents, defined as
$T^{-1} {\rm Re}\,\log \la_{\pm}$,
are given by $-3\mu C \e= -3\mu \g$.
This shows that for primary resonances of the system (\ref{eq:2.1}),
at least for initial conditions close enough to the attractors, 
convergence to the attractors has rate $1/\g$. In principle the analysis
can be extended  to any resonance, by writing $\g=C\e^{m}$ and
going up to order $m$, for suitable $m$ depending on the resonance
($m=n(\om)$ for the resonance with frequency $\om$; see Section \ref{sec:2}):
the contributions to the Lyapunov exponents due
to the Hamiltonian components of the vector fields cancel out
and the leading part of the remaining part turns out proportional to $-\g$.

In the case of  the linearly increasing friction (\ref{eq:3.2}) one expects that
the Lyapunov exponent be still proportional to $-\g$. If the friction increases very slowly,
one may reason as it were nearly constant over long time intervals,
that is time intervals covering several periods,
by approximating $\g(t)$ with a piecewise constant function.
For each of such interval $\g$ can be considered as constant and one can
reason as above. When passing from an interval to another,
the value of the initial phase $t_{0}$ of the attractor slightly changes.
The Lyapunov exponent is then expected to behave proportionally to
\begin{equation}
- \lim_{\tau \to \infty} \frac{1}{\tau} \int_{0}^{\tau} \der t \, \g(t) , \qquad
\int_{0}^{\tau} \der t \, \g(t) =
\int_{0}^{T_{0}} \der t \, \g(t) + \g_{0} \left( \tau - T_{0} \right) =
\frac{\Delta}{2} + \left( \g_{0} \tau - \Delta \right) ,
\nonumber
\end{equation}
where we have used that $T_{0}=\Delta/\g_{0}$. Therefore again
the rate of exponential convergence to the attractor is proportional to $1/\g_{0}$ and
after the time $T_{0}$ the distance to the attractor has already decreased
by a factor $\exp (-c_{0}\Delta)=\exp (-c_{0}\g_{0} T_{0})$,
for some positive constant $c_{0}$, and hence like in the case
of constant friction, possibly with a different constant $c_{0}$.

\zerarcounters
\section{Parameters $\boldsymbol{\e}$ and $\boldsymbol{\g}$ for the spin-orbit model}
\label{app:b}

The spin-orbit model has been extensively used in the literature to study
the behaviour of regular satellites \cite{GP,MD} --- it does not apply  to irregular
satellites, which are very distant from the planet and follow an
inclined, highly eccentric and often retrograde orbit.
The equations of motion are given by
\begin{equation} \label{eq:b.1}
\ddot \theta + \e \, G(\theta,t) = 0 ,
\end{equation}
with $G$ as in \eqref{eq:4.2}. Here time has been rescaled $t\to \om_{T}t$,
where $\om_{T}$ is the mean angular velocity of the satellite along its elliptic orbit
(cf. Table \ref{tab:b.1}), so that the orbital period (`year') of the satellite becomes 1.
Then the 1:1 resonance is $\dot\theta \approx 1$.

In a system satellite-primary there can be several types of friction: 
for instance the friction between the satellite layers of
different composition, say one liquid and one solid (core-mantle friction),
or the friction due to the tides (tidal friction).
One can expect that such phenomena produce a friction to be minimised in a 1:1 resonance. 
There could be also other sources of friction which we do not consider, especially those
which could modify the revolution motion of the satellite, because we are implicitly
using that it occurs on a fixed orbit. The dissipation term due to
tidal torques is of the form \cite{MacDonald,GP,Peale,CL1}
\begin{equation} \label{eq:b.2}
- \g \, \big( \Omega(e)\dot \theta - N(e) \big) ,
\end{equation}
where $\Omega(e)$ and $N(e)$ are two constants depending on $e$.
Since both $\Omega(e)=1+O(e^2)$ and $N(e)=1+O(e^2)$, for small values of $e$
we can approximate \eqref{eq:b.1}  as in \eqref{eq:4.1}.

\begin{table}[ht]
\centering
\caption{Values of $\om_{T}$ (angular velocity), $M$ (satellite mass), $M_{0}$ (primary mass),
$R$ (satellite radius) and $\rho$ (mean distance between satellite and primary) 
for the systems considered in Section \ref{sec:4}. CGS units are used.}
\begin{center}
\setlength\tabcolsep{5pt}
\vskip-.5truecm
\vrule
\begin{tabular}{lllllllllllllllllll}
\hline\noalign{\smallskip}\hline
\null & \vrule\vrule\vrule &
Primary & \vrule\vrule\vrule & Satellite & \vrule\vrule\vrule & $\om_{T}$ & \vrule & $M$ & 
\vrule & $M_{0}$ & \vrule & $R$ & \vrule & $\rho$ \\
\hline\noalign{\smallskip}\hline
E-M & \vrule\vrule\vrule & 
Earth & \vrule\vrule\vrule & Moon & \vrule\vrule\vrule
& $\st 2.66 \times 10^{-6}$ & \vrule & $\st7.35 \times 10^{25} $ & \vrule
& $\st5.97 \times 10^{27} $ & \vrule & $\st1.74 \times 10^{8} $ & \vrule & $\st3.84 \times 10^{10} $ \\
\hline
S-M & \vrule\vrule\vrule &  
Sun & \vrule\vrule\vrule & Mercury & \vrule\vrule\vrule 
& $\st 8.27 \times 10^{-7}$ & \vrule & $\st3.30 \times 10^{26} $ & \vrule
& $\st1.99 \times 10^{33} $ & \vrule & $\st2.44 \times 10^{8} $ & \vrule & $\st5.79 \times 10^{12} $ \\
\hline
J-G & \vrule\vrule\vrule &  
Jupiter & \vrule\vrule\vrule & Ganymede & \vrule\vrule\vrule 
& $\st 1.02 \times 10^{-5}$ & \vrule & $\st1.48 \times 10^{26} $ & \vrule
& $\st1.90 \times 10^{30} $ & \vrule & $\st2.63 \times 10^{8} $ & \vrule & $\st1.07 \times 10^{11} $ \\
\hline
J-I & \vrule\vrule\vrule &  
Jupiter & \vrule\vrule\vrule & Io & \vrule\vrule\vrule 
& $\st 4.11 \times 10^{-5}$ & \vrule & $\st8.93 \times 10^{25} $ & \vrule
& $\st1.90 \times 10^{30} $ & \vrule & $\st1.82 \times 10^{8} $ & \vrule & $\st4.22 \times 10^{10} $ \\
\hline
S-E & \vrule\vrule\vrule &  
Saturn & \vrule\vrule\vrule & Enceladus & \vrule\vrule\vrule 
& $\st 5.31 \times 10^{-5}$ & \vrule & $\st1.08 \times 10^{23} $ & \vrule
& $\st5.68 \times 10^{29} $ & \vrule & $\st2.52 \times 10^{7} $ & \vrule & $\st2.38 \times 10^{10} $ \\
\hline
S-D & \vrule\vrule\vrule &  
Saturn & \vrule\vrule\vrule & Dione & \vrule\vrule\vrule 
& $\st 2.66 \times 10^{-5}$ & \vrule & $\st1.09 \times 10^{24} $ & \vrule
& $\st5.68 \times 10^{29} $ & \vrule & $\st5.62 \times 10^{7} $ & \vrule & $\st3.77 \times 10^{10} $ \\
\hline
\end{tabular}
\hspace{-0.1cm}\vrule
\vskip-.5truecm
\end{center}
\label{tab:b.1}
\end{table}

A comparison with the literature \cite{GP,CL1} gives
\begin{equation} \label{eq:b.3} 
\e = \frac{3}{2} \frac{I_{y}-I_{x}}{I_{z}} \approx \frac{3h}{2R} , \qquad \g =
\frac{3k_{2}}{\xi Q} \Big( \frac{R}{\rho} \Big)^{3} \Big( \frac{M_{0}}{M} \Big) , 
\end{equation}
where $I_{x},I_{y},I_{z}$ are the moments of inertia of the satellite,
$h$ is the maximal equatorial deformation
(tide excursion), $R$ and $M$ are the mean radius and the mass of the satellite, $M_{0}$
is the mass of the primary, $\rho$ is the mean distance between
satellite and primary and $k_{2},\xi,Q$ are constants, known
respectively as the potential Love number, the structure constant and the quality factor.
For instance, in the case of the Moon one has $k_{2} \approx 0.02$, $\xi \approx 0.4$
and $Q \approx 30$ \cite{KMWL,Zhang}, which gives $3k_{2}/\xi Q \approx 0.005$
and hence $\g \approx 3.75 \times 10^{-8}$ (approximately $3.15 \times 10^{-6}$ years$^{-1}$).
In the case of S-M, the constants are usually (somehow arbitrarily) set equal to the values
$k_{2} \approx 0.4$, $\xi \approx 0.3333$ and $Q \approx 50$ \cite{GP,GS,SSWC,HSHVDS},
which gives $3k_{2}/\xi Q \approx 0.072$. The corresponding value of
the damping coefficient is $\g \approx 3.24 \times 10^{-8}$, a value very close to Mercury's.
Expressed in years$^{-1}$ this becomes approximately $8.46 \times 10^{-7}$
(the revolution period of Mercury is $2\pi/\om_{T} \approx 7.60 \times 10^{6}$ s $\approx 0.24$ year).
For lack of astronomical data we set $3k_{2}/\xi Q=10^{-1}$
for all other primary-satellite systems considered in Section \ref{sec:4}: the corresponding 
values of $\g$ as obtained from \eqref{eq:b.3} are given in Table \ref{tab:b.2}.
Of course such values only provide a rough guide.

\begin{table}[ht]
\centering
\caption{Values of $T$ (orbital period) and $\g$
for the systems considered in Section \ref{sec:4},
with $3k_{2}/\xi Q =0.1$ for the systems with Jupiter and Saturn as primary.
In the third column $\g$ is computed by
using $T$  as time unit, whereas the fourth column
gives the value of the damping coefficient expressed in years$^{-1}$.}
\begin{center}
\setlength\tabcolsep{5pt}
\vskip-.5truecm
\vrule
\begin{tabular}{llllllllllllll}
\hline\noalign{\smallskip}\hline
\null & \vrule\vrule\vrule &
Primary & \vrule\vrule\vrule & Satellite & \vrule\vrule\vrule
& $T$ &  \vrule & $T$ (years) & \vrule & $\g$ &
\vrule & $\g$ (years$^{-1}$) \\
\hline\noalign{\smallskip}\hline
E-M & \vrule\vrule\vrule & 
Earth & \vrule\vrule\vrule & Moon & \vrule\vrule\vrule
& $\st2.36 \times 10^{6}$ & \vrule & $\st7.48 \times 10^{-2}$ & \vrule
& $\st3.75 \times 10^{-8}$ & \vrule & $\st3.15 \times 10^{-6} $ \\
\hline
S-M & \vrule\vrule\vrule &  
Sun & \vrule\vrule\vrule & Mercury & \vrule\vrule\vrule 
& $\st7.60 \times 10^{6}$ & \vrule & $\st2.41 \times 10^{-1}$ & \vrule
& $\st3.24 \times 10^{-8}$ & \vrule & $\st8.46 \times 10^{-7} $ \\
\hline
J-G & \vrule\vrule\vrule &  
Jupiter & \vrule\vrule\vrule & Ganymede & \vrule\vrule\vrule 
& $\st6.18 \times 10^{5}$ & \vrule & $\st1.96 \times 10^{-2}$ & \vrule
& $\st1.91 \times 10^{-5}$ & \vrule & $\st2.48 \times 10^{-2} $ \\
\hline
J-I & \vrule\vrule\vrule &  
Jupiter & \vrule\vrule\vrule & Io & \vrule\vrule\vrule 
& $\st1.53 \times 10^{5}$ & \vrule & $\st4.84 \times 10^{-3}$  & \vrule
& $\st1.71 \times 10^{-4}$ & \vrule & $\st5.51 \times 10^{-2} $ \\
\hline
S-E & \vrule\vrule\vrule &  
Saturn & \vrule\vrule\vrule & Enceladus & \vrule\vrule\vrule 
& $\st1.18 \times 10^{5}$ & \vrule & $\st3.75 \times 10^{-3} $ & \vrule
& $\st6.26 \times 10^{-4}$ & \vrule & $\st1.05$ \\
\hline
S-D & \vrule\vrule\vrule &  
Saturn & \vrule\vrule\vrule & Dione & \vrule\vrule\vrule 
& $\st2.36 \times 10^{5}$ & \vrule & $\st7.49 \times 10^{-3}$  & \vrule
& $\st1.71 \times 10^{-4}$ & \vrule & $\st1.44 \times 10^{-1} $ \\
\hline
\end{tabular}
\hspace{-0.1cm}\vrule
\vskip-.5truecm
\end{center}
\label{tab:b.2}
\end{table}

To obtain the value of $\e$ one can use the formula
\begin{equation} \label{eq:b.4}
h = \frac{3}{2} h_{2} R \, \Big( \frac{R}{\rho} \Big)^{3} \Big( \frac{M_{0}}{M} \Big) , 
\end{equation}
for the equatorial deformation; see Appendix \ref{app:c}. Here $h_{2}$ is the tidal Love number
($h_{2} \approx 2k_{2}$ \cite{Zhang}), while the other constants are as defined after \eqref{eq:b.3}. 
If we are interested only in orders of magnitude, we can express
the equatorial deformation according to \eqref{eq:b.4} with $h_{2}=1$
for the systems with Jupiter or Saturn as primary.
Then, by inserting the values of $R,\rho,M_{0},M$ listed in Table \ref{tab:b.1}
into \eqref{eq:b.4} and using \eqref{eq:b.3} to compute $\e$, we obtain
the values in Table \ref{tab:4.1}. A comparison between
\eqref{eq:b.3} and \eqref{eq:b.4} gives $\g = C\, \e$, with $C \approx 0.05$
(taking the values of the constants $k_{2}$, $\xi$, $Q$ and $h_{2}$ in
the literature gives $C\approx 0.04$ for E-M and $C \approx 0.055$ for S-M).

Note that the values so obtained for $\e$ are lower than
those usually assumed in the literature: compare for instance 
the values $\e=2.3 \times 10^{-4}$ for E-M \cite{WND} and
$\e=1.5 \times 10^{-4}$ for S-M \cite{ACELT,CL1,CL2,CC,HSHVDS} 
with the corresponding values $\e=6.75 \times 10^{-7}$ 
and $\e=8.11\times 10^{-7}$ in Table \ref{tab:4.1}.
However, as discussed in Section \ref{sec:5}, if one does not insist
at looking only at the present structure of the satellite,
then all its evolution plays a relevant role.
So, one has to take into account that in the past,
when the satellite was more fluid, because of the lower
value of viscosity, not only the friction was smaller,
but also the deformation was bigger and hence the coupling $\e$ was larger;
see also the comments in Section \ref{sec:6}.

\zerarcounters
\section{The equatorial deformation}
\label{app:c}

Consider a homogeneous celestial body $S$ of mean radius $R$
coated by an ocean of depth $h>0$, not too small. Let $P$ be the
centre of attraction.  Denote by $M$ and $M_0$
the respective masses and assume that the motion of the two celestial
bodies about their centre of mass be circular uniform. Let $\rho$ be the
distance between the two celestial bodies $S$ and $P$, with $\rho \gg R\gg h$.  Assume,
for simplicity, that $S$ rotates about an axis orthogonal to the plane
of the orbit and that the ocean density is negligible with respect to
the core assumed to be rigid. The discussion below is essentially taken from \cite{MD}.

The distance $\rho_{C}$ of the centre of mass $C$ from the centre
of $S$ is such that $\rho_{C}(M_0+M)=M_{0}\rho$. Moreover, if
$\om_{T}$ denotes the angular velocity of revolution of the two celestial bodies
and $\kappa$ is the universal gravitational constant, one has $\om_{T}^2\rho^{3}=\kappa
(M+M_0)$ by Kepler's third law. Let $n$
be a unit vector out of the surface of $S$ and note that, imagining the observer
standing on the frame of reference rotating around $C$ with angular
velocity $\om$, so that the axis from $P$ to $S$ has a fixed unit vector $i$,
the potential (gravitational plus centrifugal) energy in the point
along the direction $n$ at distance $r$ from the centre of $S$
has density $V=V_{S}+V_{P}+V_{\rm cf}$, where
\begin{equation} \label{eq:c.1}
V_{S} = -\kappa \frac{M}r , \qquad 
V_{P} = -\kappa\frac{M_0}{(\rho^2+r^2-2 \rho r \cos\psi)^{1/2}}, \qquad
V_{\rm cf}=\frac12\om^2\,(\rho_C^2+ r^2-2 \rho_{C} r \cos\psi) ,
\end{equation}
if $\cos\psi  := i\cdot n$. Expanding $V$ in powers of $r/\rho$ one finds
\begin{equation} \label{eq:c.2}
V_{P} + V_{\rm cf} =
-\kappa\frac{M}\rho \left( \frac{r}\rho \right)^2 
\left( \frac32 \frac{M_0}{M} \cos^{2}\psi + \frac12 \right)  +
\hbox{const.} ,
\end{equation}
because the linear terms cancel out in virtue of Kepler's third law.
Therefore the equation of the equipotential surface is
\begin{equation} \label{eq:c.3}
\frac\rho{r}+ \left( \frac{r}\rho \right)^2 \left(
\frac32 \frac{M_0}{M} c^2 +\frac12 \right) = \hbox{const.} 
\end{equation}
If one writes $r=R(1+\de(\psi))$, with $\de(\psi)=\de_0+\de P_{2}(\cos\psi)$,
where $P_{2}(z)=(3z^2-1)/2$ is the second Legendre polynomial \cite{GR}
and $\de_{0}$ is such that the volume of the body $P$ is the same
as the volume of a sphere of radius $R$, then \eqref{eq:b.3} gives
\begin{equation} \label{eq:c.4}
\de(\psi) = \de P_{2}(\cos\psi) , \qquad \de_{0}=0 , \qquad 
\de = \left( \frac{R}{\rho} \right)^{3} \frac{M_{0}}M  .
\end{equation}

If the core of the satellite is rigid but the ocean density $\sigma_{\rm o}$ is not
negligible, e.g. it is equal to the core density $\sigma_{\rm c}$, then one has to take
into account that the tide will modify the potential $V_{S}$ at the site of
coordinates $r,\psi$. We make the Ansatz that the equipotential surface
is still described as $r=R (1+\de P_{2}(\cos\psi))$, possibly for a
different constant $\de$. Then the density $V_{S}$ will be
\begin{equation} \label{eq:c.5}
- \frac{3\kappa M}{4\pi R^3} \int_0^{\pi} \sin\al \,{\rm d}\al
\int_0^{2\pi} {\rm d}\varphi
\int_{0}^{R(1+\de P_{2}(\cos\al))}  \!\!\!\!\!\!\!\!\!\!\!\!\!
\frac{\rho^{2} {\rm d}\rho}{(r^{2} +
\rho^{2} -2 r\rho (\cos\psi \cos\al+\sin\psi \sin\al \cos\varphi))^{1/2}} .
\end{equation}
By expanding the integrand into Legendre polynomials and using the
orthogonality properties of the polynomials one finds
(see \cite{MD}, \S 4.3 for details)
\begin{equation} \label{eq:c.6}
V_{S} = \begin{cases}
- {\displaystyle \kappa M \left( \frac{3R^2 - r^2}{2R^3} + \frac35
\frac{r^2}{R^3} \de P_{2}(\cos\psi) \right) }  , & r<R, \\
- {\displaystyle \kappa M \left( \frac{1}{r} + \frac35
\frac{R^2}{r^3} \de P_{2}(\cos\psi) \right) } , & r>R.
\end{cases}
\end{equation}
By expanding \eqref{eq:c.6} in powers of $r/\rho$ and summing the leading
orders to \eqref{eq:c.2}, one finds that the equation of the equipotential surface becomes
\begin{equation} \label{eq:c.7}
\frac{\kappa M}{\rho} \left( \frac\rho{R} \left( 1- \frac25 \de P_{2}(\cos\psi) \right) + 
\left( \frac{r}\rho \right)^2 \frac{M_0}{M} P_2(\cos\psi) \right) = \hbox{const.} ,
\end{equation}
up to higher order corrections in $R/\rho$.
Hence if we look for the constant potential surface we find
\begin{equation} \label{eq:c.8}
\de(\psi) = \de P_{2}(\cos\psi) , \qquad \de = \frac52 \left( \frac{R}{\rho} \right)^{3} \frac{M_{0}}M  ,
\end{equation}
which replaces the previous \eqref{eq:c.4}.
The tidal deformation at the surface of the ocean, using
the notations common in celestial mechanics, can be written
as $h_2 \zeta P_2(\cos\psi)$, so that the maximal tidal excursion is
\begin{equation} \label{eq:c.9}
h = \frac32 h_2 \zeta, \qquad \zeta=R\left(\frac{R}\rho \right)^3\frac{M_0}{M}, \qquad h_2=\frac52 .
\end{equation}
The number $h_2$ is called the tidal Love number. More generally
$h_2$ depends on the detailed structure of the satellite, so far
supposed to have uniform density $\sigma_{\rm o}=\sigma_{\rm c}$.
If on the contrary $\sigma_{\rm o} \neq \sigma_{\rm c}$, then, denoting by
$r_{\rm c}$ the shape of the core boundary (while $r$ is the 
shape of the external ocean surface), one can make once more
the Ansatz that the deformations be such that $r=R(1+\de P_{2}(\cos\psi))$ and
$r_{\rm c}=R_{c}(1+ \de' P_{2}(\cos\psi))$, where $R_{\rm c}$ is
the mean radius of the core and $\de,\de'$ are two constants
to be determined by imposing that the ocean surface is equipotential
and balancing the forces acting on the core boundary. The latter
can be performed by considering the pressures acting on the core boundary
due to the elastic forces within the core and the loaded terms caused
by the ocean and core tide. This leads to two relations involving
$\de,\de'$ (see \cite{MD}, \S 4.4)
\begin{subequations} \label{eq:c.10}
\begin{align}
\frac{\z_{\rm c}}{R_{\rm c}} & =
\left[ \frac25 \frac{\sigma_{\rm o}}{\sigma_{\rm c}} + \left( \frac{R_{\rm c}}{R} \right)^3
\left( 1 - \frac{\sigma_{\rm o}}{\sigma_{\rm c}} \right) \right] \de' -
\frac35 \left( \frac{R_{\rm c}}{R} \right)^5 \left( 1 -
 \frac{\sigma_{\rm o}}{\sigma_{\rm c}} \right) \de ,
\label{eq:c.10a} \\
\de & = \frac{1}{\tilde \mu}
\left( 1 - \frac{\s_{\rm o}}{\sigma_{\rm c}} \right)
\left[ \frac52 \frac{\zeta_{\rm c}}{R_{\rm c}} - \de +
\frac32 \frac{\sigma_{\rm o}}{\sigma_{\rm c}} \left( \de' - \de \right) \right]  ,
\label{eq:c.10b}
\end{align}
\end{subequations}
\vskip-.1truecm
\noindent
where $\zeta_{\rm c}=(M/M_{\rm c})\zeta$ (with $M_{\rm c}$ being the mass of the core) and 
the effective rigidity $\tilde\m$ is a dimensionless quantity proportional to the rigidity 
of the core. For instance, if $\sigma_{\rm o} \ll \sigma_{\rm  c}$, we can approximate
\begin{equation} \nonumber
R_{\rm c} \de' \approx \frac52 \frac{\zeta_{\rm c}}{1+\tilde\mu} , \qquad
R\de \approx \zeta_{\rm c} \left[ \left( \frac{R}{R_{\rm c}} \right)^4 +
\frac32 \left( \frac{R_{\rm c}}{R} \right) \frac1{1+\tilde\mu} \right] .
\end{equation}
In particular in the limit of high rigidity ($\mu \gg 1$)
then $R_{\rm c}\de' \approx 0$ and $R\de \approx (R/R_{\rm c})^4\zeta_{\rm c}=\zeta$
(in agreement with \eqref{eq:b.4}), so that the core deformation
becomes very small, that is the core is essentially undeformed.

\zerarcounters
\section{Threshold values for the spin-orbit model}
\label{app:d}

We reason as done for the driven cubic oscillator in \cite{BBDGG}. We consider \eqref{eq:4.1}
with $\g=\e \, C$ and write it as the first order differential equation
\begin{equation} \label{eq:d.1}
\begin{cases}
\dot \theta = y , \\
\dot y = -\e \, G(\theta,t)-\e \, C \left( y - 1 \right) .
\end{cases}
\end{equation}
Then we look for a solution $z(t)=(\theta(t),y(t))$ in the form of a
power series in $\e$, that is $z(t)=z^{(0)}(t)+\e\,z^{(1)}(t) + \e^{2}z^{(2)}(t)+\ldots$,
where $z^{(0)}(t)=(\theta_{0}+\om t, \om)$, with $\om=p/q$, and
$z^{(k)}(t)=(\theta^{(k)}(t),y^{(k)}(t))$ to be determined by imposing
that $z(t)$ be periodic in $t$ with period $2\pi q$.

A first order analysis gives
\begin{equation} \label{eq:d.2}
\begin{cases}
\dot \theta^{(1)} = y^{(1)} , & \cr
\dot y^{(1)} = - G(\theta_{0}+\om t , t) - C \left( \om - 1 \right) .
\end{cases}
\end{equation}
By introducing the Wronskian matrix
\begin{equation} \label{eq:d.3}
W(t) = \left( \begin{matrix} 1 & \om t \\ 0 & \om  \end{matrix} \right) ,
\end{equation}
we can write $z^{(1)}(t)$ as
\begin{equation} \label{eq:d.4}
\left( \begin{matrix} \theta^{(1)}(t) \\ y^{(1)}(t) \end{matrix} \right) =
W(t) \left( \begin{matrix} \bar \theta^{(1)} \\ \bar  y^{(1)} \end{matrix} \right) +
W(t) \int_{0}^{t} {\rm d}\tau \, W^{-1}(\tau)  \left( \begin{matrix} 0 \\ 
-G(\theta_{0}+\om \tau, \tau) - C \left( \om  -1 \right)  \end{matrix} \right) ,
\end{equation}
with $(\bar\theta^{(1)},\bar y^{(1)})$ to be fixed. Then we obtain
\begin{equation} \label{eq:d.5}
\theta^{(1)}(t) = \bar \theta^{(1)} + \bar y^{(1)} \om t -
\int_{0}^{t} {\rm d}\tau \int_{0}^{\tau} {\rm d}\tau' 
\left[ G(\theta_{0}+\om \tau', \tau') + C \left( \om - 1 \right) \right] ,
\end{equation}
whereas $y^{(1)}(t)=\dot \theta^{(1)}(t)$. For \eqref{eq:d.4} to be periodic
we have to require first of all that
\begin{equation} \label{eq:d.6}
M(\theta_{0}) : = \frac{1}{2\pi q} \int_{0}^{2\pi q} {\rm d} t
\left[ G(\theta_{0}+\om t, t) + C \left( \om - 1 \right) \right] = 0 ,
\end{equation}
then fix $\bar y^{(1)}$ in such a way that
\begin{equation} \label{eq:d.7}
\bar y^{(1)} \om = \frac{1}{2\pi q} \int_{0}^{2\pi q} {\rm d}t
\int_{0}^{t} {\rm d}\tau \left[ G(\theta_{0}+\om \tau, \tau) + C \left( \om - 1 \right) \right] ,
\end{equation}
while $\bar\theta^{(1)}$ will be fixed to second order by requiring that also $\theta^{(2)}(t)$ be periodic.

Using that $G(\theta,t)=\partial_{\theta}g(\theta,t)$, with
$g(\theta,t)$ given by (\ref{eq:4.3}), and
inserting \eqref{eq:d.6} into \eqref{eq:d.5} leads to
\begin{equation} \label{eq:d.8}
\frac{1}{2\pi q} \sum_{k\in\ZZZ} 2 a_{k} \int_{0}^{2\pi q} {\rm d}t \,
\sin(2\theta_{0}+2\om t-kt) = C \left( \om -1 \right) 
\end{equation}
and hence
\begin{equation} \label{eq:d.9}
2 a_{k(p)} \sin 2\theta_{0} = C \Big( \frac{p}{q} - 1 \Big) ,
\qquad k(p) = \frac{2p}{q} .
\end{equation}
Since $a_{k} \neq 0$ only for $k=-3,\ldots,7$, $k\neq 0$, as \eqref{eq:4.2} implies,
$\om=p/q$ is either integer (and $\om\in\{-1,1,2,3\}$) or half-integer
(and $\om\in\{-3/2,-1/2,1/2,3/2,5/2,7/2\}$). If we confine ourselves to positive $\om$,
we see that \eqref{eq:d.9} fixes $\theta_{0}$ provided
\begin{equation} \label{eq:d.10}
\left| C \right| < C_{0}(p/q) : = \frac{2a_{k(p)} q}{|p-q|} .
\end{equation}
In particular \eqref{eq:d.10} is always satisfied for $\om=1$, so that
the 1:1 resonance always exists, while for the other values of $\om$ we obtain \eqref{eq:4.4}.

\zerarcounters
\section{Time evolution of friction for the spin-orbit model}
\label{app:e}

The mechanism of capture into resonance has been studied by
several authors starting with the theory of capture into the 3:2
resonance of Mercury \cite{GP,BJ,N}. Usually the
friction is considered either periodic or just not depending on time.
Here we regard the friction as not periodic in time and given by \eqref{eq:2.2},
that is starting from a initial very small value, then slowly increasing
in order of magnitude, until the satellite has completely solidified: such a situation
seems possible in the formation of a satellite or planet.
At the beginning, the satellite can be considered in a fluid state; however
the dissipation due to tidal effects becomes more and more sensible due to the
cooling and the resulting increase of viscosity and, eventually,
it settles at the final present value: the time over which the entire
process takes place is called the solidification time and will be denoted by $T_{S}$.

We stress that, in the model we are considering, we assume that the
satellite has first stabilised in its orbit around the primary and then
modifies its spinning velocity. Of course the exact evolution of the
satellite dynamics is still debated and no theory is universally agreed upon. 
In what follows, we shall ignore the model-dependent details of such an evolution.
So, for instance, if we accept that in some stage of the history of Mercury
large quantities of its mantle material have been removed
\cite{Woolfson2,BAHW}, for the purposes of our argument it is
not important whether this occurred before Mercury attained its final
orbital motion or after that event. 

Suppose that friction is essentially due to tidal effects on an originally entirely fluid 
fast rotating satellite (that is with rotation frequency $\dot\theta$ at least a few times larger than
the present-day orbital frequency $\om_{T}$) evolving toward a solid body. 
Assuming that the dynamics is described by \eqref{eq:4.1} since the beginning
could appear contradictory with the fact that originally the satellite was essentially fluid,
since \eqref{eq:4.1} deals with a rigid body with given moments of inertia.
On the other hand, as we shall argue, friction becomes really effective only when the
viscosity has attained high values and, when this occurs, the shape
of the satellite can be considered close to its final state. One could
also imagine to study the deformations of the satellite during its
evolution, but of course this would make the analysis much more
complicated; see for instance \cite{BH}, where asymptotic stability of the 1:1 resonance 
is obtained for a deformable body with high rigidity (see also \cite{Shatina}). 
In other words, using the spin-orbit model is justified
except possibly for the very early stage of the satellite evolution
(which are not interested in).

In the early history of the solar system one can assume that
the viscosity is rather low, not too different from that of the water,
which equals $10^{-2}$ poises (CGS units). The solidification
process is very complicated, and although it has been extensively
studied in the literature, especially in the case of the Earth
and the Moon (for the obvious reason that it is much easier
to compare the results obtained by theoretical methods and
numerical simulations with the experimental data), still there
are many unsolved issues. See for instance \cite{Shearer}
for a review of recent results on the Moon.
Usually one assumes that originally all satellites and planets were completely or
almost completely molten, according to the so-called magma ocean hypothesis;
in the case of the Moon see \cite{Warren} and especially \cite{Shearer} and the references
quoted therein (see also \cite{Turcotte} in the case of Earth). Most of the celestial body
crystallised, from the bottom to the top, following a sequence
determined by chemical composition and pressure \cite{STN,KPMH,ET} leaving only a
layer of very fluid magma close to the surface of the satellite.
The outer liquid part eventually disappeared, through a further solidification process from the
core-mantle boundary to the surface, and it is irrelevant
in any case to the damping, because of its very low viscosity.

Timing of formation of satellites and for the cooling of the magma
is mostly deduced from the study of rock samples: in the case of the Moon
of course this is much easier \cite{Shearer,Brandon,DBYJ,NTPGRM},
while the results are far from being conclusive in the case of Mercury
(however see, for instance, \cite{Solomon,BET}).
In any case,  it is not unlikely that the solidification
time of Mercury is of the same order of magnitude of that of the Moon;
we can also mention that the theory has been proposed that
the two celestial bodies may have had similar origin \cite{Woolfson1,Woolfson2}.
There is strong evidence that the solidification time $T_{S}$ for the
Moon is of order of $10^{8}$ years \cite{Brandon,NTPGRM}.
Thus, for the reasons given above, we take this as the order of
magnitude of the solidification time of all the satellites considered in Section \ref{sec:4}.

For most of the solidification process the ocean
magma is maintained above its solidus \cite{Abe,ET}.
So it is reasonable that the viscosity has increased as an
effect of the cooling of the satellite: eventually most of the satellite
is almost solidified and its viscosity is very large, say of the order of that of
the mantle, which is almost solid (hence about $10^{24}$ poises,
much higher than the viscosity of the terrestrial magma which ranges
from $10^{4}$ and $10^{12}$ poises \cite{PMC,KZ}).
Of course the viscosity strongly depends on temperature,
which in turn decreases in time during the process of cooling of the satellite.
One could look at the literature for profiles of temperatures versus
time \cite{Turcotte,FRS} or for the dependence of viscosity on temperature in
fluids \cite{Seeton}. A detailed discussion
on viscosity evolution during the solidification of the satellite
stands as a very hard problem, and in fact such an issue
is not widely studied \cite{PS,ET,Shearer}; see however \cite{RCTCSG}, where the
despinning of Saturn's satellite Iapetus is studied and
an Arrhenius law is assumed to describe the temperature
dependence of the viscosity, and \cite{RTCCM}, where the despinning of Mercury
is related to the thermal evolution of the planet. 

As far as the satellite can be considered essentially fluid,
the damping coefficient $\g$ has to be formed with the
following physical quantities: the viscosity $\eta$ of the magmatic
fluid constituting the satellite, the mass $M$ of the satellite, the
equatorial deformation $2h$ due to the tides
and the angular velocity $\om_T$. Starting from the rescaled equation \eqref{eq:4.1},
one expects $\g\propto \h R^2 h/I_{z}$, with $I_{z} \propto MR^{2}$;
then, by rescaling time $t \to t/\om_T$,  an appropriate choice for $\g$ turns out to be
$\g \sim \eta h/M\om_T$.
The evolution of magma from the initial melt passes through the formation of cumulate rocks
and fractional crystallisation, leading to the reduction of the outer
liquid layer with low viscosity and the accretion of the internal,
highly viscous core; for a more detailed discussion see for instance
\cite{Shearer} in the case of the Moon. Therefore, when the
solidification process attained a high stage of advancement,
a different model for the friction has to be taken. If the satellite
is essentially solid, with a thin external fluid layer (ocean), one can
neglect the fluid part, because of its very low viscosity, and
concentrate on the solid core. One can use the analysis in Appendix \ref{app:b}, 
in the limit case $\sigma_{\rm o} \ll \sigma_{\rm c}$ and very high rigidity $\mu$,
so that the core deformation in \eqref{eq:c.10b} is very small. Also in
this case, one can express $\g$ in terms of the involved physical quantities
and, again on the basis of dimensional arguments, set $\g \sim \mu h/M\om_{T}^{2}$.
Therefore, to sum up, friction increases with viscosity up to a certain value.
Once such a value is reached, at which the satellite can be considered
essentially a solid with very high rigidity.

The despinning time $T_{D}$ represents the time which
the satellite needs to be really attracted into a resonance. Hence
(see Appendix \ref{app:a}) $T_{D}=O(\g^{-1})$, 
where $\g$ is the quantity appearing in \eqref{eq:4.1}.
Of course, we would want that the solidification time be larger or at least
comparable with despinning time. This is the case if
one can assume $T_{S} \approx 10^{8}$ years and $T_{D} \approx 1/\g$,
with $\g$ as in Table \ref{tab:b.2} (which yields $T_{D} \approx 10^{6}$ years).

\zerarcounters
\section{Numerical details}
\label{app:f}

The numerical results have been found by running a variety of computer
programs which implement different algorithms. The main algorithms used were
(i) a standard Runge-Kutta integrator with automatic step-size control~\cite{numrec};
(ii) the Bulirsch-Stoer algorithm~\cite{numrec}, which extrapolates the step-size to zero
and is suitable for high-accuracy computation; and (iii) a numerical
implementation of the Frobenius method. Of these, (ii)
is the slowest, but serves to confirm the results obtained from the other methods.
Most of the data in this paper came from (i) and (iii), of which
(iii) is the fastest. This is because the use of series often enables large time steps
to be taken. On the assumption that the solution to the differential equation around
a point $t = t_0$ can be expressed as a power series in $t$, we obtain a (nonlinear) 
recursion formula for the coefficients in this series. In principle, as many terms as 
desired can be computed, and in practice about 25 worked well. For a
given initial condition, this series enables us to compute $x$ and
$\dot{x}$ for $|t-t_0|< R$, for some $R$ that depends on the desired
accuracy, the initial conditions, and $t_0$. The size of the step, $t - t_0$, is
chosen as the one that makes
the absolute value of the right-hand side of the differential equation,
which should of course be zero, smaller than some tolerance: $10^{-12}$ was used here.
Typical step sizes ranged between about 0.29 and 10.0, with 0.7 being
chosen about 50\% of the time. Even the smallest step size is many times larger
than that which would be used in a Runge-Kutta implementation.

The initial conditions are chosen randomly and are uniformly distributed inside the square $\calQ$.
Since we seek detailed estimates of the relative areas of the basins of attraction,
at least up to the first decimal digit whenever possible, we have to take many
initial conditions: certainly 1000 initial conditions, as in \cite{CL1,CC} is not enough.
For the data to be reliable to the first decimal place, with a 95\%
confidence level, we found that $1\,000\,000$ initial conditions need to be
taken inside the square $\calQ$ defined in Section \ref{sec:2}.
A statistically rigorous justification for this confidence level,
given the number of initial conditions, can be found in~\cite[chapter 9]{wmm}.
Thus, for the relative areas in Table \ref{tab:2.3}, we used up to $1\,000\,000$
initial conditions except for smallest values of $\g$
($500\,000$ initial condition for $\g=0.00005$,
$150\,000$ initial condition for $\g=0.00001$ and
$50\,000$ initial conditions for $\g=0.000005$).
Analogously we considered $1\,000\,000$ initial conditions
for $\e=0.5$ (Table \ref{tab:5.1}) and $500\,000$
initial conditions for $\e=0.01$ (Table \ref{tab:5.2}).
Also in Section \ref{sec:3} we considered as many initial data as possible:
$1\,000\,000$ initial conditions for larger $\g$ ($\g=0.015$),
$500\,000$ for $\g=0.005$ and $250\,000$ for the smallest value $\g=0.0005$.
As a general rule, the smaller $\g$ is, the longer the integration time and hence,
by taking $\g$ smaller, we also need to diminish the number
of initial conditions, in order that the computation time does not become
prohibitively long. However, the error on the
relative areas becomes larger. That said, the general scenario
described in Section \ref{sec:3} seems clear and well supported by the numerics.

As already noted in \cite{RMG},
for $\g$ very small, the basins of attraction are spread out over the entire phase space
and become very sparse. Thus, if one wants to detect which
basin a given initial condition belongs to, very high numerical precision
is needed.

Finally, the integration time $T_{\rm int}$ must be chosen in such a way that all trajectories
reach the attractor (within a reasonably fixed accuracy).
For instance, one can take $T_{\rm int}=N/\g_{0}$, with $N=20$.
Therefore, in order to investigate the dynamics for very small
values of the damping coefficient, $T_{\rm int}$ has to be very large. 

The conclusion is that we have to follow the trajectories of a large number of
initial conditions, for very long times and with very high accuracy.
Of course, this is at odds with
obtaining results within a a reasonable time, so that we need to reach a compromise.
This has led us to the choice described above.


\end{document}